\documentclass[12pt]{article}

\usepackage{amssymb,amsmath,latexsym}

\usepackage{endnotes}
\let\footnote=\endnote
\let\enotesize=\normalsize
\def\notesname{Endnotes}%
\def\makeenmark{$^{\theenmark}$}
\def\enoteformat{\rightskip0pt\leftskip0pt\parindent=1.75em
	\leavevmode\llap{\theenmark.\enskip}}

\usepackage{natbib}
\bibpunct[, ]{(}{)}{,}{a}{}{,}%
\def\bibfont{\small}%
\def\bibsep{\smallskipamount}%
\def\bibhang{24pt}%

\usepackage{booktabs}
\usepackage{hyperref}
\usepackage{verbatim}
\usepackage{multirow}

\usepackage{graphicx}
\usepackage{epstopdf}
\usepackage{algorithm, algpseudocode}
\usepackage{capt-of}

\usepackage{longtable}
\usepackage{tabu}

\usepackage{subfig}

\usepackage[title,titletoc,toc]{appendix}

\usepackage{sectsty}
\sectionfont{\fontsize{12}{15}\selectfont}
\usepackage[font=footnotesize]{caption}
\usepackage{authblk}

\newtheorem{theorem}{Theorem}
\newtheorem{definition}{Definition}[section]
\newtheorem{proposition}{Proposition}

\newtheorem{lemma}{Lemma}[section]

\newenvironment{pack_enum}{
\begin{enumerate}
  \setlength{\itemsep}{1pt}
  \setlength{\parskip}{0pt}
  \setlength{\parsep}{0pt}}{\end{enumerate}
}

\newenvironment{pack_item}{
\begin{itemize}
  \setlength{\itemsep}{1pt}
  \setlength{\parskip}{0pt}
  \setlength{\parsep}{0pt}}{\end{itemize}
}

\newcommand{\BeginPf}{\noindent {\em Proof: }} 
\newcommand{\EndPf}{\hfill $\square$ }     

\newcommand{\R}{\mathbb{R}}

\newcommand{\bD}{\mathbb{D}}
\newcommand{\bE}{\mathbb{E}}

\newcommand{\bH}{\mathbb{H}}

\newcommand{\cE}{\mathcal{E}}
\newcommand{\cF}{\mathcal{F}}

\newcommand{\cL}{\mathcal{L}}
\newcommand{\cK}{\mathcal{K}}

\newcommand{\cN}{\mathcal{N}}

\newcommand{\fB}{\mathfrak{B}}
\newcommand{\fC}{\mathfrak{C}}

\newcommand{\fU}{\mathfrak{U}}

\newcommand{\fZ}{\mathfrak{Z}}

\begin{document}
\title{Reference-Based Almost Stochastic Dominance Rules with Application in Risk-Averse Optimization}

\author[1]{Jian Hu\thanks{jianhu@umich.edu}}
\author[1]{Gevorg Stepanyan\thanks{gstepany@umich.edu}}

\affil[1]{Department of Industrial and Manufacturing Systems Engineering, University of Michigan - Dearborn, Dearborn, MI 48128}


\date{\today}
\maketitle
\abstract{Stochastic dominance is a preference relation of uncertain prospect defined over a class of utility functions. While this utility class represents basic properties of risk aversion, it includes some extreme utility functions rarely characterizing a rational decision maker's preference. In this paper we introduce reference-based almost stochastic dominance (RSD) rules which well balance the general representation of risk aversion and the individualization of the decision maker's risk preference. The key idea is that, in the general utility class, we construct a neighborhood of the decision maker's individual utility function, and represent a preference relation over this neighborhood. The RSD rules reveal the maximum dominance level quantifying the decision maker's robust preference between alternative choices. We also propose RSD constrained stochastic optimization model and develop an approximation algorithm based on Bernstein polynomials. This model is illustrated on a portfolio optimization problem. }

\section{Introduction}
\label{sec: Introduction}

Stochastic dominance is a risk averse stochastic ordering approach (\cite{Hadar1969}, \cite{Hanoch1969}, \cite{Bawa(1975)StochasticDominance}, \cite{Shaked1994}, \cite{Muller2002}). Consider two random variables $ X, Y \in (\Omega, \cF, P; \Theta) $ with the support $\Theta := [\underline \theta, \ \overline \theta]$. $X$ stochastically dominates $Y$ in the $m^{th}$ order if $\bE[u(X)] \ge \bE[u(Y)]$  for any utility function $u$ satisfying $(-1)^{k-1} u^{(k)}(x) \ge 0$, $k = 1, \dots, m$, for all $ x \in \Theta$ (\cite{Levy:06}, \cite{Brockett:1987}). This utility class represents basic properties of risk aversion. For example, the decision maker prefers more to less if $u'(x) \ge 0$, is risk averse if $u''(x) \le 0$, and becomes prudent if $u'''(x) \ge 0$ on $\Theta$. On the one hand, the use of stochastic dominance to compare alternatives avoids assessing the decision maker's specific utility function which characterizes her risk attitude. Fully learning individual risk attitude is restricted with cognitive difficulty and incomplete information (\cite{Karmarkar:1978}, \cite{Weber(1987)UtiliySetDecision}). On the other hand, stochastic dominance based preference is often unnecessarily over-conservative (\cite{Leshno:2002}, \cite{LizyayevRuszczynski(2011):AlmostStochasticDominance}, \cite{Hu.etal(2013)StochasticallyWeightedDominance}). We consider the following example analogous to one given by \cite{Levy:06}. Suppose that Hannah wants to invest in one of the following lottery tickets priced at \$1:
\begin{pack_item}
	\item $\bar X_{1}$: yielding $ \$ 0 $ with the probability of $ 1\% $ and \$2 with the probability of $ 99\%$,
	\item $\bar Y_{1}$: yielding \$0.01 (1 penny) with certainty.
\end{pack_item}
Although $\bar X_{1}$ is much more attractive than $\bar Y_1$, stochastic dominance does not support the preference of $\bar X_{1}$ over $\bar Y_1$. In this case, the support $\Theta = [0, 2]$, where $\underline \theta = 0$ means that Hannah loses all investment and $\overline \theta = 2$ shows that Hannah gains 100\%
profit. Let Hannah's utilities at \$0 and \$2 be 0 and 1, respectively. It can be seen that, for the utility function $\tilde u(x) := \sqrt[1000]{x/2} $, we have that $\bE[\tilde u(\bar X_{1})] \le \bE[\tilde u (\bar Y_1)]$. Since $\tilde u$ is infinitely differentiable  in $ \Theta $ ($\tilde u'(0) = \infty$ is allowed)  and has the derivatives with alternating signs regarding the degrees of the derivatives, $\tilde u $ belongs to the utility class of any order stochastic dominance. The risk characterization specified by $\tilde u$ overvalues very small gains but completely neglects the possible difference in large gains ($\tilde u(x)$ has a very stiff increase for a small $x$, while being flat for a large $x$). Such behavior is quite unnatural; however, stochastic dominance requires that the preference should hold for all suitable utility functions including $\tilde u$.

In the paper we propose a novel concept of reference-based almost stochastic dominance (RSD). This concept provides a natural way to relax stochastic dominance by combining general and individual characterizations of risk aversion. The key idea is that, in the general utility classes, we construct a neighborhood of the decision maker's individual utility function, and specify preference rules over this neighborhood. This relaxation is motivated by the utility theory based interpretation of stochastic dominance. Differently, \cite{Leshno:2002}, \cite{LizyayevRuszczynski(2011):AlmostStochasticDominance}, and \cite{Hu.etal(2013)StochasticallyWeightedDominance} relaxed the distributionally defined forms of the first and second order stochastic dominance. Working on the CVaR based interpretation of the second order stochastic dominance, \cite{NoyanRudolf(2013)CVaRStochasticOrder} proposed the relaxed stochastic ordering which requires the CVaR of the preferred random variable to be larger over a shrunk set of confidence levels.  \cite{ArmbrusterDelage(2011)UtilitySet} and \cite{HuMehrotra(2015)RobustRAUtility} developed robust expected utility maximization models which also consider individualizing the set of utility functions to meet the decision maker's risk attitude.  \cite{ArmbrusterDelage(2011)UtilitySet} used the paired game method where the decision maker's risk attitude is partially characterized by her pairwise preference of designed lotteries. \cite{HuMehrotra(2015)RobustRAUtility}  specified boundary conditions on utility and marginal utility functions using parametric utility assessments, and construct auxiliary conditions using both standard and paired game methods. 

We address optimization problems with RSD constraints. In the literature, \cite{denrus:03,Dentcheva2004} first introduced stochastic dominance constrained optimization problems, which pursue expected profit while hedging risk by choosing options preferable to a random benchmark.  Since the last decade, optimization models using stochastic dominance have been the subject of theoretical considerations and practical applications in areas such as business, finance, energy and transportation (e.g. \cite{Karoui:06, RomanDarbyMitra(2006)SD_Portfolio, denrus:06, Dentcheva:07, luedtke:08, LeanMcAleerbWong(2010)SD_OilMarket, DrapkinSchultz(2010)FirstOrderSD_Alg, Hu.etal2009SampleAverageApproximation,  NieWuHomem(2012)SD_Transportation, SunXu(2014)SD_SAA, HaskellJain(2015)SD_MarkovDecisionProcess}). 

In addition, our model uses the concept of functional robustness. We specify 
a nonparametric shape-preserving utility neighborhood. This specification is suitable for classical nonparametric standard gamble methods and paired gamble methods such as preference comparison, probability equivalence, value equivalence, and certainty equivalence (\cite{Farquhar:1984}, \cite{Wakker:1996}, and reference therein). The functionally robust optimization was first proposed by \cite{HuLiMehrotra(2015)FunctionallyRobustNewsvendor} in a newsvendor problem for the unknown mathematical form of the price-demand function, which is different from
traditional robust approaches requiring the knowledge of the functional form (e.g. \cite{Ben-Tal1998,ElGhaouiLebret(1997)LeastSquaresDataUncertaintyRO,Bertsimas.etl(2004)RobustLinearNormUncertainty,Scarf:1958,Shapiro2004,DelageYe(2010)DistributionallyRO,Bertsimas.etl(2010)MinMaxDistribution}). To specify the uncertainty set of the price-demand function,  \cite{HuLiMehrotra(2015)FunctionallyRobustNewsvendor} consider an error allowance for the least-squares fitting at discrete data points. We generalize their approach to introduce an $\cL_2$-norm based perturbation tolerance around the decision maker's reference utility function. In the context of stochastic dominance, this tolerance is interpreted as the decision maker's desired dominance level.

This paper is organized as follows. In Section 2, we define RSD, and use the example of Hannah's comparing lotteries to illustrate the desired dominance level. Section 3 develops an optimization model with RSD constrains, and discusses its approximation using Bernstein polynomials. We analyze the approximation error in relation to the degree of Bernstein polynomials. In Section 4, to solve the approximation problem, we develop a cut-generation algorithm. The effect of the RSD constraint and the complexity of the algorithm are illustrated in Section 5 by using the financial portfolio optimization problem given in \cite{denrus:03}. Section 6 concludes.

\section{Reference-Based Almost Stochastic Dominance (RSD)}
\label{sec:RUBASD}
In this section we discuss the concept of RSD, which is specified as a preference relation based on the neighborhood of the decision maker's reference utility function. Without loss of generality, assume that the reference utility function, denote by $u_{ref}$, is increasing on $\Theta$ and satisfies $u_{ref}(\underline{\theta}) = 0$ and $u_{ref}(\overline{\theta}) = 1$. 
\begin{definition}
	\label{def:ASD}
	For the reference utility function $u_{ref}$, a random variable $ X \in (\Omega, \cF, P; \Theta) $ is preferred to another random variable $ Y \in (\Omega, \cF, P; \Theta) $ in the $ m^{th}$ order reference-based almost stochastic dominance (RSD) for a given $\epsilon \in [0, 1]$ (written as $ X \succeq^{\epsilon}_{(m)} Y $ w.r.t. $ u_{ref} $), if
	$$ 
	\bE[u(X)] \geq \bE[u(Y)], 
	$$
	for any utility function $u$ satisfying 
	\begin{pack_enum}
		\renewcommand{\theenumi}{{(A\arabic{enumi})}}
		\item \label{con:rasd_1} for any $x \in \Theta$, $(-1)^{i-1} u^{(i)}(x) \ge 0$, $i = 1, \dots, m$.
		\item \label{con:rasd_2} $ u(\underline{\theta}) = u_{ref}(\underline{\theta}) = 0 $, and $ u(\overline{\theta}) = u_{ref}(\overline{\theta}) = 1$,
		\item \label{con:rasd_3} $ 
		\|u-u_{ref}\|_{\cL_2}
		\leq 
		\epsilon$ ($\|f\|_{\cL_2} := \left(\int_{\Theta} f(t)^2 \mu(dt) \right)^{1/2}$ for a given nonnegative measure $\mu$ with $\mu (\Theta) = 1$),
		\item for any $x \in \Theta$, \label{con:rasd_4} $ u(x) \leq \frac{M}{1-\epsilon} u_{ref}(x) $, for a fixed $M > 1$ (if $\epsilon = 1$, then $u(x) < \infty$).
	\end{pack_enum}
\end{definition}

Condition \ref{con:rasd_1} represents the utility class required by the $m^{th}$ order stochastic dominance to describe the basic properties of risk aversion. Condition \ref{con:rasd_2} is the normalization of utility functions. Condition \ref{con:rasd_3} generates a neighborhood around reference utility function $u_{ref}$. The neighborhood is a closed ball, with the radius $\epsilon$, on the $\cL_2$-normed space with the measure $\mu$ on $\Theta$.  Note that conditions \ref{con:rasd_1} and \ref{con:rasd_2} ensure that, if $\epsilon = 1$, condition \ref{con:rasd_3} becomes redundant, and RSD in this case represents stochastic dominance. To further interpret $\epsilon$, we need to define the maximum dominance level as follows.  
\begin{definition}
	\label{def:lvlRSD}
	For two random variables $X, Y \in (\Omega, \cF, P; \Theta)$ satisfying $\bE[u_{ref}(X)] \ge \bE[u_{ref}(Y)]$, the maximum level of $X$ almost dominating $Y$ in the $m^{th}$ order w.r.t $u_{ref}$ is 
	\begin{align*}
		\cE^{(m)} (X, Y; u_{ref}) := \sup \; \{ \epsilon \in [0, 1] :  X \succeq^{\epsilon}_{(m)} Y \ \text{w.r.t.} \ u_{ref} \}.
	\end{align*}
\end{definition}

The maximum dominance level quantifies how $X$ is robustly preferred to $Y$. Note that $X$ stochastically dominates $Y$ in the $m^{th}$ order if and only if $\cE^{(m)} (X, Y; u_{ref}) = 1$. By Definition \ref{def:lvlRSD}, we interpret $\epsilon$ in condition \ref{con:rasd_3}  as the decision maker's desired dominance level with which she can assert $X$ is sufficiently preferred to $Y$ in the sense that the ambiguity and inconsistency in the elicitation of $u_{ref}$ is not very sensitive. We now illustrate the maximum dominance level using the case of Hannah's purchasing lottery tickets $\bar X_{1}$ and $\bar Y_1$. It has been shown that stochastic dominance is unable to reveal the preference of $\bar X_{1} $ over $\bar Y_1$, for conservatively taking unreasonable utility functions (e.g. $\tilde u$ given in the introduction) into consideration. Now suppose that Hannah's risk preference is approximately characterized as $\bar{u}_{ref}(x) = \sqrt{x/2}$, which we use as the reference utility function.  The maximum dominance level $\cE^{(2)}(\bar X_{1}, \bar Y_1; \bar u_{ref}) = 0.398$. Hence, with the given $\epsilon \le 0.398$, we have that $\bar X_1 \succeq^{\epsilon}_{(2)} \bar Y_1, ~w.r.t. \bar{u}_{ref}$, and condition \ref{con:rasd_3} excludes utility functions that Hannah is unwilling to choose (e.g. $\tilde u$). To understand the statement that the maximum dominance level quantifies the preference level in the sense of robustness, we consider the non-purchase option, denoted by $\bar Y_2$, for which Hannah has no gain but never takes a risk. The rational decision maker does not purchase $\bar Y_1$ in which the investment is lost for sure. This undoubted fact is supported by $\cE^{(2)}(\bar Y_2, \bar Y_1; \bar u_{ref}) = 1$, which indicates $\bar Y_2$ stochastically dominates $\bar Y_1$. So we can claim that Hannah absolutely prefers $\bar Y_2$ to $\bar Y_1$.  In contrast, although $\bar X_1$ could be better than $\bar Y_2$ ($\bar u_{ref}(\bar X_1) = 0.99 > \bar u_{ref}(\bar Y_2) = 0.707$), we cannot say that $\bar X_1$ is absolutely preferred to $\bar Y_1$ according to stochastic dominance rules. Indeed, there is still a theoretical possibility that $\bar Y_1$ is preferred to $\bar X_1$. It can only be said that the preference of $\bar X_1$ over $\bar Y_1$ is highly robust. We also consider an additional \$1 lottery ticket as
\begin{pack_item}
	\item  $\bar X_2$: yielding $ \$ 0 $ with the probability of $ 10\% $ and \$2 with the probability of $ 90\%$,
\end{pack_item}
While $\bar X_1$ stochastically dominates $\bar X_2$, $\bar X_2$ is still more attractive then $\bar Y_1$. This preference is visible also from  $\cE^{(2)}(\bar X_2, \bar Y_1; \bar u_{ref}) = 0.069 $. Comparing the two pairs $(\bar X_1, \bar Y_1)$ and  $(\bar X_2, \bar Y_1)$, we have that $\cE^{(2)}(\bar X_1, \bar Y_1; \bar u_{ref}) > \cE^{(2)}(\bar X_2, \bar Y_1; \bar u_{ref})$. These results show that $\bar X_1$ is more robustly preferred to $\bar Y_1$ than $\bar X_2$.

\begin{table}[htbp]
	\centering
	\caption{Dominance Levels of Lottery Tickets}
	\begin{tabular}{crrcc}
		\addlinespace
		\toprule
		Lottery   &   \multicolumn{2}{c}{Probability of Yield}  &  & Maximum Dominance Level  \\
		Ticket & \multicolumn{1}{r}{\$0} &  \multicolumn{1}{r}{\$2} &  & $\cE^{(2)}(\bar X_{i}, \bar Y_2; \bar u_{ref})$ \\
		\midrule 
		$\bar X_{1}$    &   1\% & 99\%   &  & 0.122    \\
		$\bar X_{2}$    &   10\% & 90\%  &  & 0.069   \\
		$\bar X_{3}$    &   25\% & 85\%  &  & 0.045   \\
		$\bar X_{4}$    &   20\% & 80\%  &  & 0.024   \\
		$\bar X_{5}$    &   25\% & 75\%  &  & 0.007  \\
		\bottomrule
	\end{tabular}%
	\label{tab:lotteries}%
\end{table}%

We next illustrate the estimation of Hannah's desired dominance level by comparing the non-purchase option $\bar Y_2$ and lottery tickets $\bar X_i$ priced at \$1,  for $i = 1, \dots, 5$, given in Table \ref{tab:lotteries}. Note that, since $\bar X_i$ stochastically dominates $\bar X_{i+1}$, the preference of $\bar X_i$ is monotonously weakened such that the maximum dominance level decreases.  Hannah is requested to choose the lottery tickets from the list which she is not reluctant to purchase. Since $\bar Y_2$ is a non-yielding but risk-free option, Hannah's choice indicates the level of her insistence on risky investment. Suppose that she picks first three lotteries unhesitatingly, but feels it difficult to make a decision on $\bar X_4$. Then her desired dominance level should be $0.045$, which is the maximum dominance level for $\bar X_3$. In other words, she would like to invest in the lottery ticket that almost dominates the non-purchase option with the maximum level no less than $0.045$. Hannah's decision is that  $\bar X_1$, $\bar X_2$, and $\bar X_3$ are sufficiently preferred to $\bar Y_2$ in 2nd order RSD, while $\bar X_4$ and $\bar X_5$ may be indifferent.

The measure $\mu$ in condition \ref{con:rasd_3} is used to quantify the relative importance of different subintervals of $\Theta$. A special case is that $\mu$ is a discrete measure and condition \ref{con:rasd_3} is the weighted least squares fitting criterion. Actually, nonparametric utility assessments can only generate finitely many utility value points, and then use a piecewise linear curve to link all these points. We may specify a perturbation set based on those discrete points instead of the piecewise linear curve. Let $(x_i, u_{ref}(x_i))$, $i = 1, \dots, I$, be the reference utility points. The measure $\mu$ should be assigned on $x_i$ and condition \ref{con:rasd_3} is thus represented as
$$ \left(\sum_{i=1}^I \mu(x_{i}) (u(x_{i}) - u_{ref}(x_{i}))^2  \right)^{1/2} \leq \epsilon.$$
Now consider the set of utility functions described by conditions \ref{con:rasd_1} - \ref{con:rasd_3}. The closure of this set contains utility functions which are discontinuous at the lower boundary point $\underline \theta$. Condition \ref{con:rasd_4} gives a point-wise upper bound of the utility set. The constant $M$ in the condition can be a very larger number.  With this weak condition, the specified set is closed in the space of continuous function defined in $\Theta$. Moreover,  condition \ref{con:rasd_4} may only exclude the characterizations of extremely irrational risk attitudes. Note that condition \ref{con:rasd_4} is redundant for stochastic dominance when $\epsilon = 1$. 

Let
\begin{align}
	\fU^m(\epsilon) := \{ u : u(x) \ \text{satisfies conditions \ref{con:rasd_1}-\ref{con:rasd_4} for given  $\epsilon$} \}.
\end{align}
By Definition \ref{def:ASD}, $ X \succeq^{\epsilon}_{(m)} Y, ~w.r.t. ~u_{ref}$ is equivalent to $\bE[u(X)] \ge \bE[u(Y)]$ for all $u \in \fU^m(\epsilon)$.   The following propositions show the relation between different orders of RSD, and claim the nonemptiness and convexity of the set $\fU(\epsilon)$.
\begin{proposition} \label{pp:orders_SD}
	If $ X \succeq^{\epsilon}_{(m)} Y, ~w.r.t. ~u_{ref}$, then $X \succeq^{\epsilon}_{(m+1)} Y, ~w.r.t. ~u_{ref}$.
\end{proposition}
\BeginPf
It follows by condition \ref{con:rasd_1}-\ref{con:rasd_4} that $\fU^{m+1}(\epsilon) \subseteq \fU^m(\epsilon)$. 
\EndPf

\begin{proposition} \label{pp:cU_convex}
	$\fU^m(\epsilon)$ is a convex set.  If $u_{ref}$ satisfies condition \ref{con:rasd_1}, then $\fU^m(\epsilon)$ is nonempty.
\end{proposition}
\BeginPf
The reference utility function $ u_{ref} $ satisfies conditions \ref{con:rasd_2}-\ref{con:rasd_4} by the construction, so if $u_{ref}$ satisfies condition \ref{con:rasd_1}, then $ u_{ref} \in \fU^m(\epsilon)$. Thus $ \fU^m(\epsilon) $ is nonempty.

We now check the convexity. For $ u_{1}, \ u_{2} \in \fU^m(\epsilon) $, let $ \tilde{u}(x) := \lambda u_{1}(x) + (1-\lambda) u_{2}(x) $ for $\lambda \in [0, 1]$ in $ \Theta $. Obviously, $ \tilde{u} $ satisfies the conditions \ref{con:rasd_1}, \ref{con:rasd_2} and \ref{con:rasd_4}. We next check condition \ref{con:rasd_3}. 
$$ \| \tilde{u} - u_{ref} \|_{\cL_2} \leq \lambda \| u_{1} - u_{ref} \|_{\cL_2} + (1 - \lambda) \| u_{2} - u_{ref} \|_{\cL_2} \leq \epsilon. $$

Thus, $ \tilde{u} \in \fU^m(\epsilon) $, and hence $ \fU^m(\epsilon) $ is convex.
\EndPf

\section{Optimization Model using the Reference-Based Almost Stochastic Dominance}
\label{sec:OMutRUBASD}
In this section we first develop a RSD constrained stochastic optimization model. For the risk-averse decision maker, utility functions should be increasing and concave. Hence, in the later statement, the reference utility function $u_{ref}(x)$ is assumed to be increasing and concave in $\Theta$, and the constraint is specified using the second or higher order RSD, i.e., the order of dominance $m \ge 2$. We next study an approximation approach to the RSD constrained model using Bernstein polynomials. Finally, we discuss the connection of the RSD constrained model with its approximation, and verify the asymptotic convergence of the approximation.

\subsection{RSD Constrained Optimization model}
\label{sec:optimizationModel}
A stochastic optimization model using the $m^{th}$ ($m \ge 2$) order RSD as a risk-averse constraint is specified as follows:
\begin{equation*}
	\label{pr:RSDP}
	\tag{RSD-P}
	\begin{aligned}
		& \underset{z}{\text{max}}
		& & f(z) \\
		& \text{subject to}
		& & X(z) \succeq^{\epsilon}_{(m)} Y, ~w.r.t. ~u_{ref}, \\
		& 
		& & z \in \fZ, \\
	\end{aligned}
\end{equation*}
where $\fZ \subseteq \R^d$ is a decision region, $f : \R^d \mapsto \R$ is an objective function, $X : \fZ \mapsto (\Omega, \cF, P; \Theta)$ represents a random outcome function of the decision, and $Y  \in (\Omega, \cF, P; \Theta)$ is a benchmark. In model \ref{pr:RSDP}, the RSD constraint requires that the random outcome $X(z)$ at a valid decision $z$ should almost dominate the random benchmark $Y$ in the $m^{th}$ order for the given dominance level $\epsilon$ with respect to the reference utility function $u_{ref}$. 

We now state notions needed in the discussion. Denote
\begin{align}
	\label{def:robustPreference}
	\pi^m (z) := \min_{u \in \fU^m(\epsilon)} \left\{ \Pi(z, u) := \bE[u(X(z))-u(Y)] \right\},
\end{align}
and
\begin{align}
	\Psi^m (\delta) := \{z \in \fZ : \pi^m(z) \ge \delta \}.
\end{align}
By Definition \ref{def:ASD}, the RSD constraint in model \ref{pr:RSDP} equals to $\pi^m (z) \geq 0$. Hence, $\Psi^m (0)$ is the feasible region of model \ref{pr:RSDP}. By Proposition \ref{pp:orders_SD}, we have the following relationship of the set $\Psi^m(\delta)$ for different $m$'s. 
\begin{proposition}
	For any $\delta \in \R$, $\Psi^m (\delta) \subseteq \Psi^{m+1} (\delta)$.
\end{proposition}

\subsection{Approximation using Bernstein Polynomials}
Model \ref{pr:RSDP} is a functionally robust optimization problem, where the RSD constraint is specified using the set $\fU^m(\epsilon)$ of nonparametric utility functions. We now discuss an approximation approach using Bernstein polynomials.

Let the vector $\phi(x) := (\phi_{0}(x), \dots, \phi_{n}(x))^{T}$, where  
$$ \phi_{j}(x) := 
\left( 
\begin{array}{c} 
n \\ j 
\end{array} 
\right) 
\left( \frac{x-\underline{\theta}}{\overline{\theta}-\underline{\theta}} \right)^{j}
\left( 1-\frac{x-\underline{\theta}}{\overline{\theta}-\underline{\theta}} \right)^{n-j}, \quad j = 0, \dots, n,$$
are the bases of Bernstein polynomial on $ \Theta $. The 
$n^{th}$ degree Bernstein polynomial is given as
$$
B_{n}(x; c) := c^{T}\phi(x), ~x \in \Theta,
$$
where $c := (c_{0}, \dots, c_{n})^{T}$ is a vector of coefficients. Note that condition \ref{con:rasd_1} requires utility functions to be $m$ times differentiable. In order to avoid trivial solutions, we require that $n \geq m$. Consider the following conditions on coefficients $c = (c_{0}, \dots, c_{n})^{T} $:
\begin{pack_enum}
	\renewcommand{\theenumi}{{(B\arabic{enumi})}}
	\item \label{con:basd_1} $ (-1)^{i-1} \Delta^i c_{j} \geq 0 $,  $ i=1 \dots m,~j = 0 \dots n-i $, where
	$$
	\Delta^i c_{j} := \sum_{k=0}^{i} (-1)^{k} 
	\left( 
	\begin{array}{c} 
	k \\ i 
	\end{array} 
	\right) c_{j+i-k}, 
	$$
	
	\item \label{con:basd_2} $ c_{0} = u_{ref}(\underline{\theta}) = 0$, and $ c_{n} = u_{ref}(\overline{\theta}) =  1$,
	
	\item \label{con:basd_3} $ c^T A c + g^{T} c + r \leq \epsilon^{2} $, where 
	\begin{align*}
		A := \int_{\Theta}
		\phi(x) 
		\phi^{T}(x)
		\mu(dx), \ 
		g := -2\int_{\Theta}
		u_{ref}(x)\phi(x)
		\mu(dx), \nonumber \ \text{and }
		r := \int_{\Theta}
		u_{ref}^{2}(x)
		\mu(dx). 
	\end{align*}
	
	\item \label{con:basd_4} $ c_{j} \leq \min \left\{1,  ~\frac{M}{1-\epsilon} u_{ref} (\underline{\theta}+j \frac{\overline{\theta}-\underline{\theta}}{n})\right\}, ~j = 0 \dots n$.
\end{pack_enum}
Denote the set of coefficient
\begin{align} \label{set:coef}
	\fC^m_{n}(\epsilon) := \{c \in \R^{n+1} : c \text{ satisfies conditions \ref{con:basd_1}-\ref{con:basd_4}}  \}.
\end{align}
The theorem below states that $\fU^m(\epsilon)$ contains the set of Bernstein polynomials with the coefficients belonging to $\fC^m_{n}(\epsilon)$. 
\begin{theorem}
	\label{thm:B_n-U}
	Let $\fB^m_n(\epsilon) := \{B_n(x; c) : c \in \fC^m_n(\epsilon) \}$.
	$\fB^m_{n}(\epsilon) \subseteq \fU^m(\epsilon)$.
\end{theorem}

\BeginPf
For any Bernstein polynomial $B_{n}(\cdot; c) \in \fB^m_n(\epsilon)$, we have $c \in \fC^m_n(\epsilon)$. We now prove that $B_n(x; c)$ satisfies conditions \ref{con:rasd_1}-\ref{con:rasd_4} in $ \Theta $, such that $B_n(x; c) \in \fU^m_n(\epsilon)$.

\ref{con:basd_1}. Theorem 7.1.3 in \cite{Philips(2003)InterpolationAndApproximationByPolynomials} shows that
$$
\frac{d^i B_{n}(x; c)}{d x^i} = \frac{n!}{(\overline{\theta} - \underline{\theta})^i(n-i)!} \sum_{j=0}^{n - i} \Delta^i c_{j} \phi_{j}(x).
$$
The Bernstein polynomial base $ \phi(x) $ is nonnegative in $ \Theta $. Hence, condition \ref{con:basd_1} ensures that  $B_{n}(x; c)$ satisfies condition \ref{con:rasd_1}.

\ref{con:basd_2}. Since $
B_{n}(\underline{\theta}; c) 
=
c_{0}
$ and $
B_{n}(\overline{\theta}; c) 
=
c_{n}
$, conditions \ref{con:rasd_2} and \ref{con:basd_2} are equivalent.

\ref{con:basd_3}. By condition \ref{con:rasd_3}, we have
\begin{align*}
	& \int_{\theta}
	(B_{n}(x; c)-u_{ref}(x))^{2} 
	\mu(dx) \\
	= \ & \int_{\theta} (c^T \phi(x) -u_{ref}(x))^{2} \mu(dx) \\
	= \ & 
	c^T 
	\left( \int_{\theta}
	\phi(x) 
	\phi^{T}(x)
	\mu(dx) \right)
	c 
	- 2 \left(
	\int_{\theta}
	u_{ref}(x)\phi(x)
	\mu(dx)
	\right)^{T}
	c 
	+ \int_{\theta}
	u_{ref}^{2}(x)
	\mu(dx) \\
	= \ & 
	c A c^{T} 
	+ g^{T} c + r \\
	\le \ & \epsilon^2.
\end{align*}
Hence,  $B_{n}(\cdot; c)$ satisfies condition \ref{con:rasd_3}.

\ref{con:basd_4}. Let 
$$
B_{n}^{ref}(x) := \sum\limits_{j=0}^{n} u_{ref}
\left(
\underline{\theta}+j\frac{\overline{\theta}-\underline{\theta}}{n}
\right) \phi_{j}(x).
$$
Since $u_{ref}$ is concave, it follows by Theorem 7.1.8 in \cite{Philips(2003)InterpolationAndApproximationByPolynomials} that 
$$
B_{n}^{ref}(x) \le u_{ref}(x), \quad \text{for all } x \in \Theta.
$$
Also, condition \ref{con:basd_4} ensures that 
$$
B_n(x; c) \le \min \left\{1, ~\frac{M}{1-\epsilon} B_{n}^{ref}(x) \right\}  \le \frac{M}{1-\epsilon} B_{n}^{ref}(x).
$$
Therefore, $B_n(x; c)$ satisfies condition \ref{con:rasd_4} on $ \Theta $. 
\EndPf

Letting 
\begin{align}
	\label{pr:pi^m_n}
	\pi^m_{n}(z) := \min_{c \in \fC^m_{n}(\epsilon)} \left\{ \Pi_n(z, c) := \bE[B_{n}(X(z), c)-B_{n}(Y, c)] \right\},
\end{align}
we now present an approximation of the RSD constraint using Bernstein polynomials as
\begin{align}
	\pi^m_{n}(z) \ge 0. \label{pr:BDP} \tag{BSD}
\end{align}
Similarly, the set $\Psi^m(\delta)$ and model \ref{pr:RSDP} are approximated as 
\begin{align}
	\Psi^m_n(\delta) := \{ z \in \fZ :  \pi^m_{n}(z) \ge \delta \},
\end{align}
and
\begin{align}
	& \max_{z \in \Psi_{n}^m (0)}	f(z).
	\label{pr:BSDP}
	\tag{BSD-P}
\end{align}
The degree $n$ of Bernstein polynomials is an important parameter of the 
approximation model \ref{pr:BSDP}. The next section will discuss the asymptotic convergence of model \ref{pr:BSDP} to model \ref{pr:RSDP} as $n$ increases to infinity. Hence, we call $n$ the degree of model \ref{pr:BSDP}.

\subsection{Relationship between Models \ref{pr:RSDP} and \ref{pr:BSDP}}
We now discuss the connection between model \ref{pr:RSDP} and its approximation \ref{pr:BSDP}. Two theorems are given below. Theorem \ref{thm:approximation} describes the relationship between the feasible regions of models \ref{pr:RSDP} and \ref{pr:BSDP}, and Theorem \ref{thm:convergence} shows the asymptotic convergence of the optimal value and the set of optimal solutions of the approximation model \ref{pr:BSDP}.   

\begin{theorem} 
	\label{thm:approximation} 
	Suppose $u_{ref} \in \fU^m(\epsilon)$.
	For any $ \delta \in (0, 1)$, choose the degree $n$ of model \ref{pr:BSDP} such that  
	\begin{align*}
		\delta  \ge \Lambda(n) := \frac{4M}{1-\epsilon} u_{ref}
		\left(
		\frac{1}{\sqrt{n}}+\underline{\theta}
		\right)
		\left(
		1+\frac{\overline{\theta}-\underline{\theta}}{n}
		\right).
	\end{align*} 
	Then
	\begin{align*}
		\Psi^m_n(\delta) \subseteq \Psi^m(0) \subseteq \Psi^m_n(0) \subseteq \Psi^m(-\delta).
	\end{align*}  
\end{theorem}

To prove theorem \ref{thm:approximation}, we need two technical lemmas. By adjusting Theorem 2.1 in \cite{Rivlin(1981)AnIntroductionToTheApproximationOfFunctions}, we obtain Lemma \ref{lem:1}, which gives the precision of Bernstein polynomial based approximation of $u \in \fU^m(\epsilon)$. The result straightforwardly follows the proof of Theorem 2.1 in \cite{Rivlin(1981)AnIntroductionToTheApproximationOfFunctions}, and  we omit the proof. Lemma \ref{lem:apprxomiation} considers the uniform bound of the approximation error. For $u \in \fU^m(\epsilon)$, we define the operators
$$
T_j u := u \left( \underline \theta + \frac{j (\overline \theta - \underline \theta)}{n}\right), \quad j = 0, \dots, n,
$$ 
and let 
$$
T u := (T_1 u, \dots, T_n u).
$$

\begin{lemma}[Theorem 1.2 in \cite{Rivlin(1981)AnIntroductionToTheApproximationOfFunctions}] \label{lem:1}
	
	For $\tilde \delta > 0$, the modulus of continuity of $u \in \fU^m(\epsilon)$ is 
	$$
	\omega(\tilde \delta) := \sup_{
		\tiny{\begin{array}{c}
			x_1, x_2 \in \Theta \\
			|x_1 - x_2| \le \tilde \delta
			\end{array}
		} } |u(x_1) - u(x_2)|.
		$$
		Then, 
		$$ \left\| u - B_n(\cdot ; T u) \right\|_\infty \le \omega \left( \frac{1}{\sqrt{n}} \right) \left( 1 + \frac{\overline \theta - \underline \theta}{n}\right).$$
	\end{lemma}

	\begin{lemma}
		\label{lem:apprxomiation}
		For $\tilde \delta > 0$, if 
		$
		\Lambda(n) \le \tilde \delta,
		$ 
		then, for $u \in \fU^m(\epsilon)$, we have 
		$$\left \| u - B_n(\cdot; T u)  \right \|_\infty \le \frac{\tilde \delta}{4}, \text{ and } T u \in \fC^m_n \left( \epsilon + \frac{\tilde \delta}{4} \right).$$
	\end{lemma}
	\BeginPf
	Condition \ref{con:rasd_1} for $m \ge 2$ ensures that $u \in \fU^m(\epsilon)$ is increasing concave on $\Theta$ such that 
	$$ 
	\omega \left( \frac{1}{\sqrt{n}} \right) 
	= u \left( \underline{\theta} + \frac{1}{\sqrt{n}} \right) - u \left( \underline{\theta} \right) 
	= u \left( \underline{\theta} + \frac{1}{\sqrt{n}} \right).
	$$
	By condition \ref{con:rasd_4}, we also have
	$$
	u \left( \underline{\theta} + \frac{1}{\sqrt{n}} \right) \le \frac{M}{1-\epsilon} u_{ref} \left( \underline{\theta} + \frac{1}{\sqrt{n}} \right).
	$$
	Hence, it follows by Lemma \ref{lem:1} that
	$$
	\left\| u - B_{n}(\cdot; Tu) \right\|_{\infty} 
	\leq 
	\frac{M}{1-\epsilon} u_{ref} \left( \underline{\theta} + \frac{1}{\sqrt{n}} \right)
	\left( 1+\frac{\overline{\theta}-\underline{\theta}}{n} \right)
	=
	\frac{\Lambda(n)}{4} \le \frac{\tilde \delta }{4}.
	$$
	
	The definition of $\fU^m(\epsilon)$ implies that $Tu$ satisfies conditions \ref{con:basd_1}, \ref{con:basd_2}, and \ref{con:basd_4}. Condition \ref{con:basd_3} can also be verified as
	\begin{align*}
		\|B_n(\cdot; Tu) - u_{ref}\|_{\cL_2} \le & \|B_n(\cdot; Tu) - u\|_{\cL_2} + \|u - u_{ref}\|_{\cL_2} \\
		\le & \left\| u - B_{n}(\cdot; Tu) \right\|_{\infty} + \epsilon \\
		\le & \epsilon + \frac{\tilde \delta  }{4}.
	\end{align*}
	\EndPf

	\BeginPf[Proof of Theorem \ref{thm:approximation}] 
	
	Let $ u^{*} $ be the optimal solution of problem \eqref{def:robustPreference}. We define
	$$ u_{\delta}(x) := \frac{\delta}{4}u_{ref}(x) + \left(1-\frac{\delta}{4} \right) u^{*}(x), ~x \in \Theta.$$
	Since condition \ref{con:rasd_2} implies that maximum difference between $u^*$ and $u_{ref}$ on $\Theta$ is no more than 1, we have that
	$$
	\|u_\delta - u^*\|_\infty = \frac{\delta}{4}  
	\| u_{ref} - u^* \|_\infty \le \frac{\delta}{4}.
	$$
	By the assumption that $u_{ref} \in \fU^m(\epsilon)$ and the convexity of $\fU^m(\epsilon)$ shown in Proposition \ref{pp:cU_convex}, we have $u_\delta \in \fU^m(\epsilon)$. Since
	\begin{align*}
		\left\|
		u_{\delta} - u_{ref}
		\right\|_{\cL_2}
		=
		\left(1-\frac{\delta}{4}\right)
		\left\|  u^{*} - u_{ref}
		\right\|_{\cL_2}
		\leq
		\left(1-\frac{\delta}{4}\right)\epsilon
	\end{align*}
	we have $u_\delta \in \fU^m \left(\left(1-\frac{\delta}{4}\right)\epsilon \right)$. By Lemma \ref{lem:apprxomiation} and the assumption that $\Lambda(n) \ge \delta \epsilon$, we have 
	$$ T u_\delta \in \fC^m_{n}(\frac{\delta\epsilon}{4} + \epsilon - \frac{\delta\epsilon}{4}) = \fC^m_{n}(\epsilon)$$ 
	and
	$$
	\| u_{\delta} - B_n(\cdot; T u_\delta) \|_{\infty} \le \frac{\delta \epsilon}{4} \le \frac{\delta}{4}.
	$$
	It follows that 
	\begin{align*}
		\left\| B_{n}(\cdot; T u_\delta) - u^{*} \right\|_{\infty}
		\leq
		\left\| B_{n}(\cdot; T u_\delta) - u_{\delta} \right\|_{\infty}
		+
		\left\| u_{\delta} + u^{*} \right\|_{\infty}
		\leq \frac{\delta}{4} + \frac{\delta}{4}
		=
		\frac{\delta}{2},
	\end{align*}
	and hence, for any $z \in \fZ$, 
	\begin{align*}
		& \ \left\| B_{n}(X(z), T u_\delta) - B_{n}(Y, T u_\delta)  - (u^{*}(X(z)) - u^{*}(Y))\right\|_{\infty} \\
		\leq & \
		\left\| B_{n}(X(z), T u_\delta) - u^*(X(z)) \right\|_{\infty}
		+
		\left\| B_{n}(Y, T u_\delta) - u^*(Y) \right\|_{\infty} \\
		\leq & \ 
		\delta. 
	\end{align*}
	We can conclude that
	$$
	\pi^m_{n}(z) \leq 
	\bE[ B_{n}(X(z), T u_\delta) - B_{n}(Y, T u_\delta)]
	\leq
	\delta + \bE[u^{*}(X(z)) - u^{*}(Y)] 
	= 
	\delta + \pi^m(z).
	$$
	On the other hand, Theorem \ref{thm:B_n-U} shows that 
	\begin{align*}
		\pi^m(z) \leq \pi^m_{n}(z). 
	\end{align*}
	By the definitions of the sets $\Psi^m$ and $\Psi^m_n$, it follows that $\Psi^m_n(\delta) \subseteq \Psi^m(0) \subseteq \Psi^m_n(0) \subseteq \Psi^m(-\delta)$.
	\EndPf

	We now discuss the asymptotic convergence of the optimal value and the set of optimal solutions of the approximation model \ref{pr:BSDP}.   
	
	\begin{theorem}
		\label{thm:convergence}
		Let $\xi^m$ and $\Xi^m$ be the optimal value and the set of optimal solutions of model \ref{pr:RSDP}, $\xi^m_n$ and $\Xi^m_n$ be the optimal value and the set of optimal solutions of model \ref{pr:BSDP} with degree $n$. Suppose that (i) $u_{ref} \in \fU^m(\epsilon)$, (ii) the set $\fZ$ is convex and compact, (iii) the random function $X(z)$ is concave in $\fZ$, (iv) there exists an interior point $\tilde z \in \fZ$ such that $\pi(\tilde z) \ge 0$, and (v) the objective function $f$ is continuous in $\Psi^m(0)$, then $\xi^m_n \to \xi^m$ and $\bD(\Xi^m_n, \Xi^m) := \max_{x \in \Xi^m_n} \min_{y \in \Xi^m} \|x - y\|_\infty \to 0$ as $n \to \infty$.  
	\end{theorem}
	\BeginPf
	For $\delta \in (0, 1)$, denote by $(\Psi^m(-\delta))^o$ the set of all the interior points of $\Psi^m(-\delta)$. Assumption (iv) ensures that $(\Psi^m(-\delta))^o$  is nonempty. Assumption (iii) implies that $\Pi(z, u)$ defined in \eqref{def:robustPreference} is concave in $\fZ$ for any $u \in \fU^m(\epsilon)$, and hence, $\pi(z)$ is concave in $\fZ$. This conclusion with assumption (ii) shows that $\Psi^m(-\delta)$ is convex. The convexity guarantees that $\Psi^m(-\delta) \in cl ((\Psi^m(-\delta))^o)$ which is the closure of $(\Psi^m(-\delta))^o$. It follows that $\bD(\Psi^m(-\delta), \Psi^m(0)) \to 0$ as $\delta \to 0$. By assumption (i) and Theorem \ref{thm:approximation}, we know that $\Psi^m(0) \subseteq \Psi^m_n(0) \subseteq \Psi^m(-\delta)$ when $n$ satisfies $\Lambda(n) \le \delta \epsilon$. It means that $\bD(\Psi^m(0), \Psi^m_n(0)) = 0$ for any $n$. On the other hand, $\bD(\Psi^m_n(0), \Psi^m(0)) \le \bD(\Psi^m(-\Lambda(n)/\epsilon), \Psi^m(0))$ for $n$ sufficiently large, and hence, 
	$\bD(\Psi^m_n(0), \Psi^m(0)) \to 0$ as $n \to \infty$. Therefore, as $n \to \infty$, the Hausdorff distance of these two sets $$\bH(\Psi^m(0), \Psi^m_n(0)) := \max \{\bD(\Psi^m(0), \Psi^m_n(0)), \ \bD(\Psi^m_n(0), \Psi^m(0))\} \to 0.$$
	
	Assumptions (ii)-(iv) ensure that $\Psi^m(0)$ is nonempty and compact. Then, $\Xi^m$ is also nonempty because of the continuity of $f$ given in assumption (v). For $z^* \in \Xi^m$, let $$z_n := \arg\min_{z \in \Xi^m_n} \|z - z^*\|_\infty.$$ 
	Since $\bH(\Psi^m(0), \Psi^m_n(0)) \to 0$, it follows that $z_n \to z^*$ as $n \to \infty$. By the continuity of $f$, 
	$$\lim\inf_{n \to \infty} \xi^m_n \ge \lim_{n \to \infty} f(z_n) = f(z^*) = \xi^m.$$ 
	Also, for any convergent sequence $\{\bar z_n\}$ with $\bar z_n \in \Xi^m_n$, the limit point $\bar z$ is in $\Xi^m$. Hence, 
	$$\lim\sup_{n \to \infty} \xi^m_n \le \xi^m.$$ 
	Suppose $\bD(\Xi^m_n, \Xi^m) \nrightarrow 0$. We can construct a convergent sequence  $\{\bar z_n\}$ with $\bar z_n \in \Xi^m_n$ such that 
	$$\min_{z \in \Xi^m} \|z - \bar z_n\|_\infty > \tau > 0.$$ 
	Let $\bar z$ be its limit point, so $\bar z \in \Psi^m(0) \setminus \Xi^m$. But 
	$$f(\bar z) = \lim_{n \to \infty} f(\bar z_n) = \lim_{n \to \infty} \xi^m_n = \xi^m,$$ 
	which is a contradiction. Thus, $\bD(\Xi^m_n, \Xi^m) \to 0$. 
	\EndPf

	\section{Cut Generation Algorithm for Model \ref{pr:BSDP} }
	\label{sec:CGAfM}
	We now develop a cut generation algorithm to solve model \ref{pr:BSDP}. This algorithm uses a sequence of coefficient cuts $c^i  \in \fC^m_n(\epsilon)$ for $i = 1, \dots, k$. We solve a sequence of the problems
	\begin{equation}
		\label{pr:r-BSDP}
		\begin{aligned}
			\max_z \ & f(z) \\
			\text{s.t.} \ & \Pi_{n}(z, c^{i}) \geq 0, \quad i = 1, \dots, k, \\
			& z \in \fZ, \\
		\end{aligned}
	\end{equation}
	which are the relaxations of model \ref{pr:BSDP} over the subset of $\Psi^m_n(0)$ consisting of the generated cuts. At the optimal solution $z^*_k$ of problem \eqref{pr:r-BSDP}, we calculate $\pi^m_n(z^*_k)$ by solving problem \eqref{pr:pi^m_n}. Denote by $c^*_k$ the optimal solution of problem \eqref{pr:pi^m_n}. If $\pi^m_n(z^*_k) \ge 0$, $z^*_k$ is the optimal solution of model \ref{pr:BSDP}. Otherwise, the constraint $\Pi_n(z, c^*_k) \ge 0$ is added to the master problem \eqref{pr:r-BSDP} as a valid cut. Algorithm \ref{alg:bernsteinRelaxedGeneralOptimizationAlgorithm} formally describes this cut generation algorithm. If $\fZ$ is convex and both $f(z)$ and $X(z)$ are concave in $\fZ$, the master problem \eqref{pr:r-BSDP} is a stochastic convex program, which can be solved via sample average approximation method  (see \cite{Shapiro:08}).  Problem \eqref{pr:pi^m_n}  is a quadratic constrained linear program (QCP)(see \cite{vandePanne(1966)QCP}, \cite{MarteinSchaible(1987)QCP}, \cite{Boyd:2004}).
	
	\algnewcommand{\Label}{\State\unskip}
	
	\begin{algorithm}[h]
		\caption{Cut-Generation Algorithm for Model \ref{pr:BSDP} \label{alg:bernsteinRelaxedGeneralOptimizationAlgorithm}}
		\begin{algorithmic}[]
			\Label \texttt{Step 1}	
			\begin{raggedright}
				\hspace*{.2in}
				\label{alg:step0}
				Choose $\gamma > 0$ and let $k = 0$. \\
				\par \end{raggedright}

			\Label \texttt{Step 2}	
			\begin{raggedright}
				\hspace*{.2in}
				\label{alg:step1}
				Find the optimal solution $z_{k}^{*}$ of problem \eqref{pr:r-BSDP}.\\
				\par 	\end{raggedright}
			
			\Label \texttt{Step 3}
			\begin{raggedright}
				\hspace*{.2in}
				\label{alg:step2}
				Calculate $\pi^m_{n}(z^*_k)$ by solving problem \eqref{pr:pi^m_n}. Let $c^*_k$ be the optimal solution. \\
				\par \end{raggedright}
			
			\Label \texttt{Step 4}
			\begin{raggedright}
				\hspace*{.2in}
				\label{alg:step3}
				If $ \pi_{n}(z_{k}^{*}) \geq -\gamma $, exit. Otherwise, let $c^{k+1} = c^*_k$ and $k = k+1$. Then go to step 1.\\
				\par \end{raggedright}
		\end{algorithmic}
	\end{algorithm}
	
	The following  theorem \ref{thm:finishingTime} shows that Algorithm \ref{alg:bernsteinRelaxedGeneralOptimizationAlgorithm} terminates in finitely many iterations. Let $$\xi^m_n (-\gamma) := \max_{z \in \Psi^m_n(-\gamma)}f(z),$$ which is a relaxation of model \ref{pr:BSDP} for $\gamma > 0$. $\xi^m_n (0)$ is the optimal value of model \ref{pr:BSDP}. 
	
	\begin{theorem}
		\label{thm:finishingTime}
		Algorithm \ref{alg:bernsteinRelaxedGeneralOptimizationAlgorithm} terminates in finitely many iterations. Let $\widetilde \xi^m_n$ be the optimal value of problem  \eqref{pr:r-BSDP} at the last iteration where Algorithm \ref{alg:bernsteinRelaxedGeneralOptimizationAlgorithm} terminates.  Then $\xi^m_n(0) \le \widetilde \xi^m_n \le \xi^m_n(-\gamma)$.
	\end{theorem}
	\BeginPf
	Let us say that Algorithm \ref{alg:bernsteinRelaxedGeneralOptimizationAlgorithm} does not terminate in $k$ iterations. Then the stopping criterion at step 4 is not satisfied, i.e., $\pi^{*}_{n}(z^{*}_{k}) = \Pi_{n}(z^{*}_{k}, c^*_k) < -\gamma$. Recall that $z^*_k$ and $c^*_k$ are optimal solution of problems \eqref{pr:r-BSDP} and \eqref{pr:pi^m_n}, respectively. 
	
	Denote the close balls on the $\cL_\infty$ space 
	$$\cN(c^i) := \left\{c \in \R^{n+1} : \|c - c^i\|_\infty \le \frac{\gamma}{2} \right\}, \quad i = 1, \dots, k.$$ 
	We claim that $c^*_k$ is not covered by any of these balls. By contradiction, suppose that there is $j \in \{1, \dots, k\}$ such that $c^*_k \in \cN(c^j)$. Then, by the H\"{o}lder's inequality and the fact that $\|\phi(x)\|_1 = 1$ for any $x \in \Theta$, we have that
	$$
	| \Pi_{n}(z^*_k, c^j) - \Pi_{n}(z^*_k, c^*_k) | 
	\leq
	\| {c^j} - {c^*_k} \|_\infty
	\| \bE[\phi(X(z_{k}^{*})) - \phi(Y)] \|_1
	\le \frac{\gamma}{2} \bE[\| \phi(X(z_{k}^{*}))\|_1 + \| \phi(Y) \|_1] = \gamma.
	$$
	Since $z^{*}_{k}$ is the optimal solution of problem \eqref{pr:r-BSDP} at the $k^{th}$ iteration, we have that $\Pi_{n}(z^{*}_{k}, c^{j}) \geq 0$. It follows
	$$
	\Pi_{n}(z^{*}_{k}, c^{*}_k) 
	\ge
	\Pi_{n}(z^{*}_{k}, c^{j}) - \gamma 
	\ge
	-\gamma,
	$$
	which contradicts to the stopping criterion. 
	
	The generated cut with respect to $c^{k+1} = c^{*}_k$ results in the addition of the open ball $\cN(c^{k+1})$ at the $(k+1)^{th}$ iteration. On the other hand, by conditions \ref{con:basd_1}, \ref{con:basd_2}, we know that $0 \le c \le 1$ for any $c \in \fC^m_n(\epsilon)$ (all the components of $c$ should be in [0, 1]). Hence, $\fC^m_n(\epsilon)$ is compact, and  $$\max_{c^1, c^2 \in \fC^m_n(\epsilon)} \|c^1 - c^2\|_\infty \le 1.$$ It means that Algorithm \ref{alg:bernsteinRelaxedGeneralOptimizationAlgorithm} at most runs $\left\lceil \left(\frac{2}{\gamma}\right)^n \right\rceil$ iterations in total.
	
	Suppose that Algorithm \ref{alg:bernsteinRelaxedGeneralOptimizationAlgorithm} terminates at the $k^{th}$ iteration. Denote by $\widetilde \Psi^m_n$ the feasible solution of problem \eqref{pr:r-BSDP} at the $k^{th}$ iteration. Then, $\widetilde \Psi^m_n$ is represented as 
	$$
	\widetilde \Psi^m_n = \{z \in \fZ : \Pi_n(z, c^i) \ge 0, \ i = 1, \dots, k, \ \Pi_n(z, c) \ge -\gamma, \text{ for all } c \in \fC^m_n(\epsilon) \}.
	$$  
	Obviously, $\Psi^m_n(0) \subseteq  \widetilde \Psi^m_n \subseteq \Psi^m_n(-\gamma)$. Then, $\xi^m_n(0) \le \widetilde \xi^m_n \le \xi^m_n(-\gamma)$.
	\EndPf

	\section{Case Study: Portfolio Investment}
	\label{sec:CS}
	
	In this section we apply the framework \ref{pr:RSDP} to the portfolio optimization problem given by \cite{denrus:03}. This problem involves $N (= 8)$ assets: (S1) U.S. three-month treasury bills, (S2) U.S. long-term government bonds, (S3) S\&P 500, (S4) Willshire 5000, (S5) NASDAQ, (S6) Lehmann Brothers corporate bond index, (S7) EAFE foreign stock index, and (S8) gold. \cite{denrus:03} use $M(=22)$ yearly returns $r_{ij} \ (i = 1, \dots, M, \  j = 1, \dots, N)$ of these assets as equally probable realizations (See Table \ref{tbl:Index_Return} in the appendix). 
	
	Using the framework \ref{pr:RSDP} to model this problem, we have
	\begin{subequations} \label{pr:portfolio}
		\begin{align}
			\max_z \ & \left(1+\frac{1}{M} \sum_{i=1}^M  \sum_{j=1}^N z_j r_{ij} \right) \label{pr:portfolio_obj}\\
			\text{s.t.} \ &  \frac{1}{M} \sum_{i=1}^M  u \left( 1 + \sum_{j=1}^N z_j r_{ij} \right) \ge \frac{1}{M} \sum_{i=1}^M  u \left( 1 + \sum_{j=1}^N z^Y_j r_{ij} \right), \quad  \text{for all } u \in \fU^m (\epsilon), \label{pr:portfolio_RSD} \\
			& \sum_{j=1}^N z_j = 1, \\
			& z_j \ge 0, \quad j = 1, \dots, N.
		\end{align}
	\end{subequations}  
	In the above problem, the objective \eqref{pr:portfolio_obj} seeks the best asset allocation to maximize the expected total wealth, while the RSD constraint \eqref{pr:portfolio_RSD} requires that this allocation should be sufficiently preferred to the benchmark $z^Y$. 
	
	We discuss model \eqref{pr:portfolio} in four cases. In case (i), we consider the 2nd order RSD constraint by setting $m = 2$ in \eqref{pr:portfolio_RSD}. The option of investing all money in S1 is used as the benchmark, i.e., $z^{Y_1} := \{1,0,\dots,0\}$. S1 is a risk-free asset, using which the RSD constraint \eqref{pr:portfolio_RSD} guarantees the investment on risky assets to reach a given level of safety. we let the support $\Theta = [0, 2]$ and choose the CRRA utility function, $u_{ref}^{1}(x) := \sqrt{\frac{x}{2}}$, as the reference, which is consistent with the example of Hannah's purchasing lottery tickets. Case (i) is default in this study, and we adapt it to the other three cases: case (ii) substitutes the CARA reference utility function $u_{ref}^{2}(x) := \frac{e^{x} - 1}{e^2 - 1}$; case (iii) uses an alternative benchmark,  $z^{Y_2} := \left\{ \frac{1}{N},\frac{1}{N},\dots,\frac{1}{N} \right\}$, which equally invests on every asset; and case (iv) discusses the 3rd order RSD constraint by letting $m = 3$. Table \ref{tab:config} summarizes the different configurations  in these cases.  
	
	\begin{table}[!htbp]
		\centering
		\caption{Configurations in the Four Studied Cases}
		\begin{tabular}{lccc}
			\addlinespace
			\toprule
			&   Reference utility & Benchmark & RSD order \\
			\midrule 
			Case (i)    & $ u_{ref}^1 $ & $ z^{Y_1} $ & 2 \\
			Case (ii)    & $ u_{ref}^2 $ & $ z^{Y_1} $ & 2 \\
			Case (iii)    & $ u_{ref}^1 $ & $ z^{Y_2} $ & 2 \\
			Case (iv)    & $ u_{ref}^1 $ & $ z^{Y_1} $ & 3 \\
			\bottomrule
		\end{tabular}%
		\label{tab:config}%
	\end{table}%
	
	Model \eqref{pr:portfolio} is approximated by the corresponding \ref{pr:BSDP}. This study first tests the computational complexity of Algorithm \ref{alg:bernsteinRelaxedGeneralOptimizationAlgorithm} with $\gamma = 10^{-10}$, and next analyzes the model performance by adjusting dominance level $\epsilon$. Finally, we discuss the application of model \eqref{pr:portfolio} to Hannah's investment. 
	
	\subsection{Computational Analysis}
	
	\begin{figure}[ht]
		\subfloat[Optimal value]{
			\label{fig:compTotalWealth}
			\includegraphics[width=3in,height=2.2in]{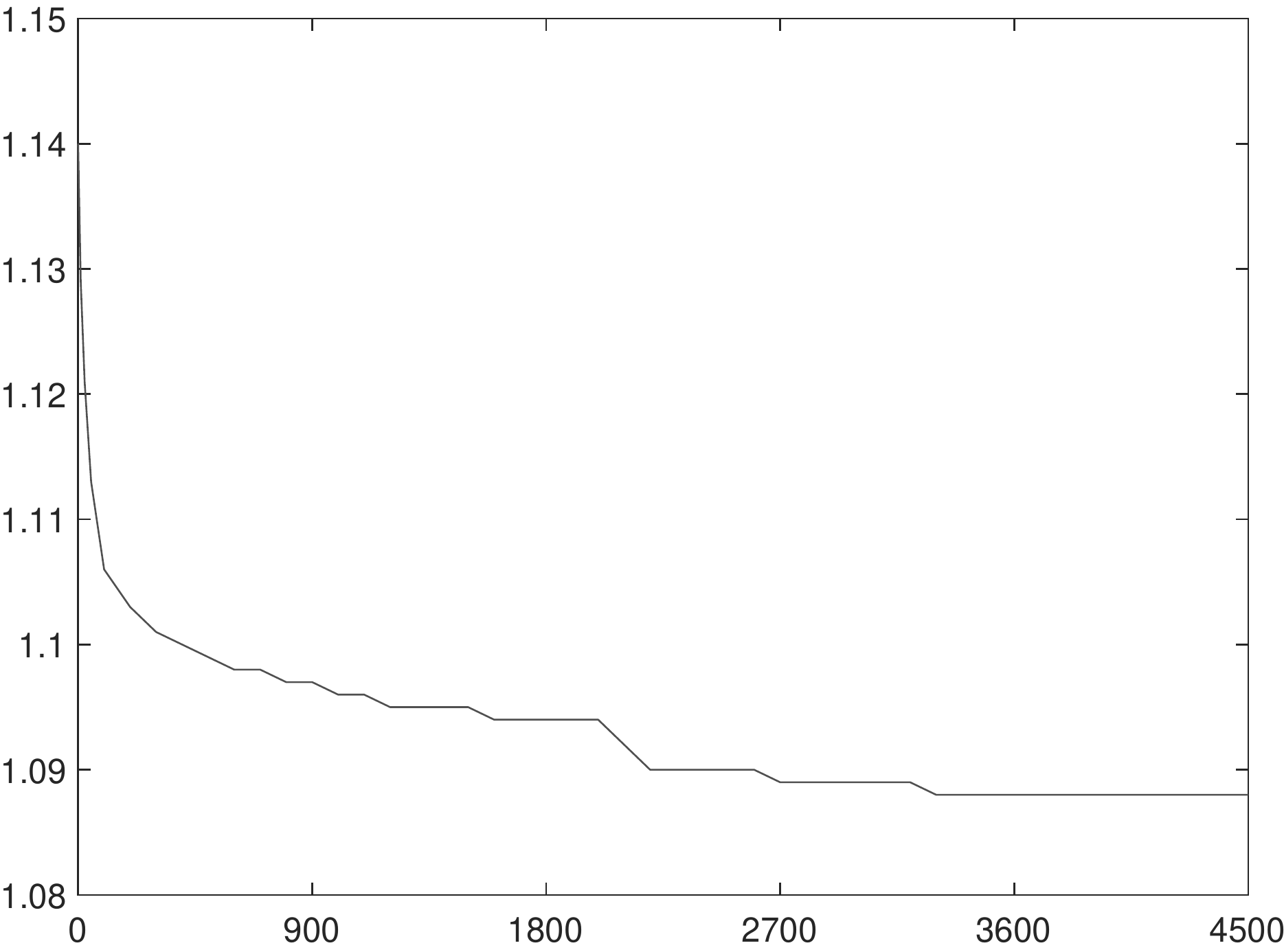}
		}
		\qquad
		\subfloat[Optimal solution]{
			\label{fig:compSol}
			\includegraphics[width=3in,height=2.2in]{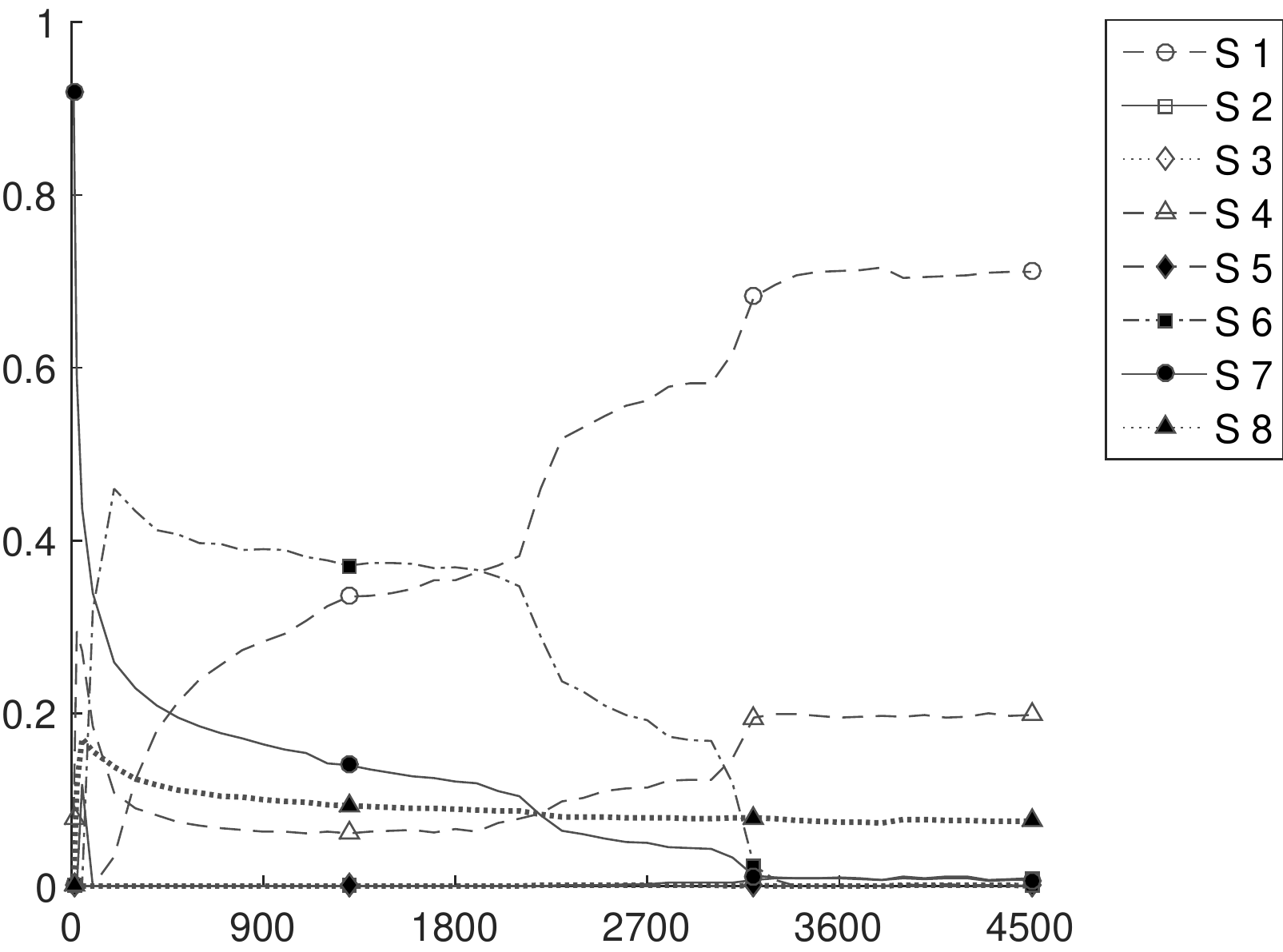}
		}
		\\
		\subfloat[Runnning time (in sec).]{
			\label{fig:compRunningTime}
			\includegraphics[width=3in,height=2.2in]{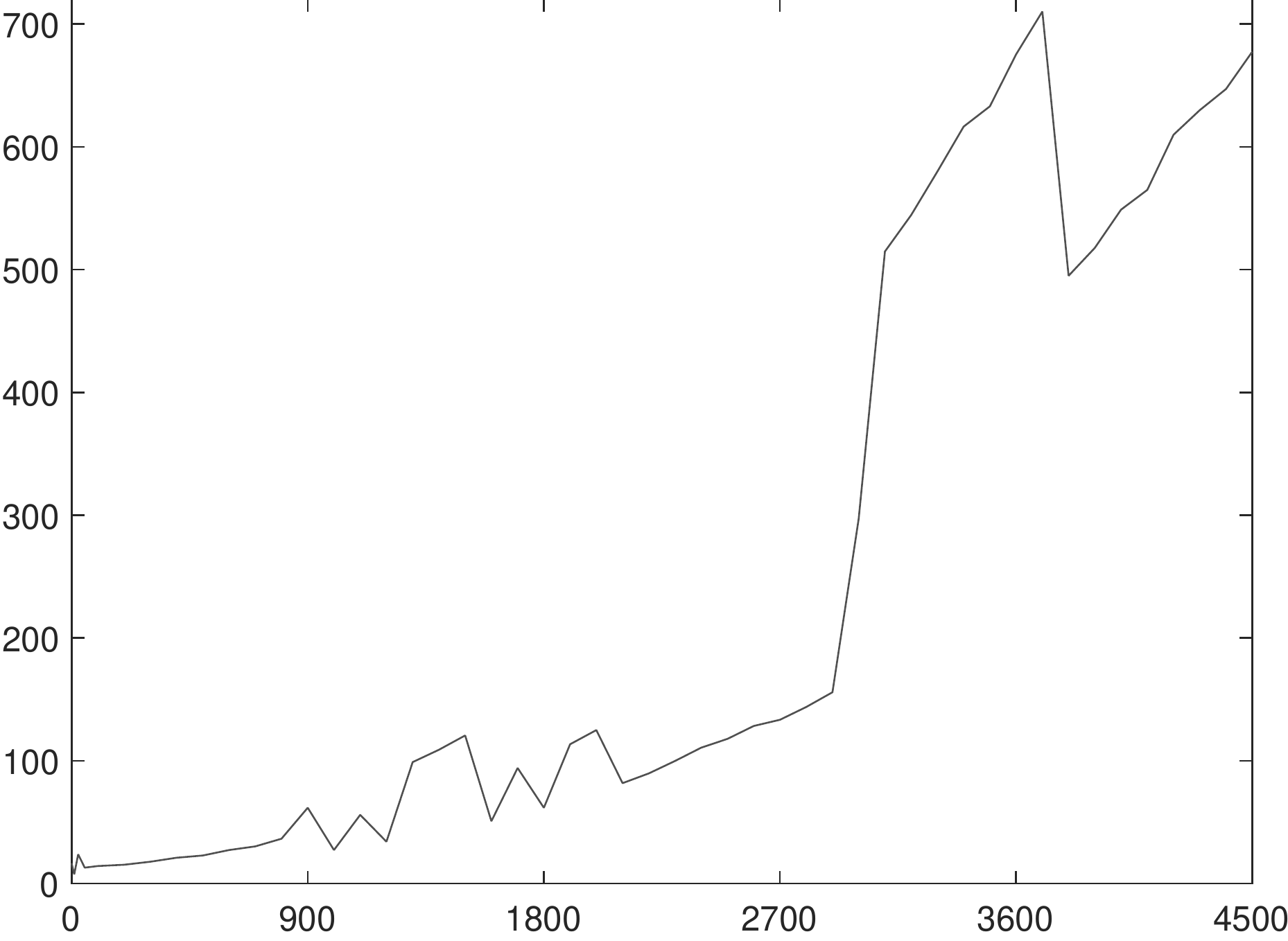}
		}
		\qquad
		\subfloat[Number of generated cuts.]{
			\label{fig:compCuts}
			\includegraphics[width=3in,height=2.2in]{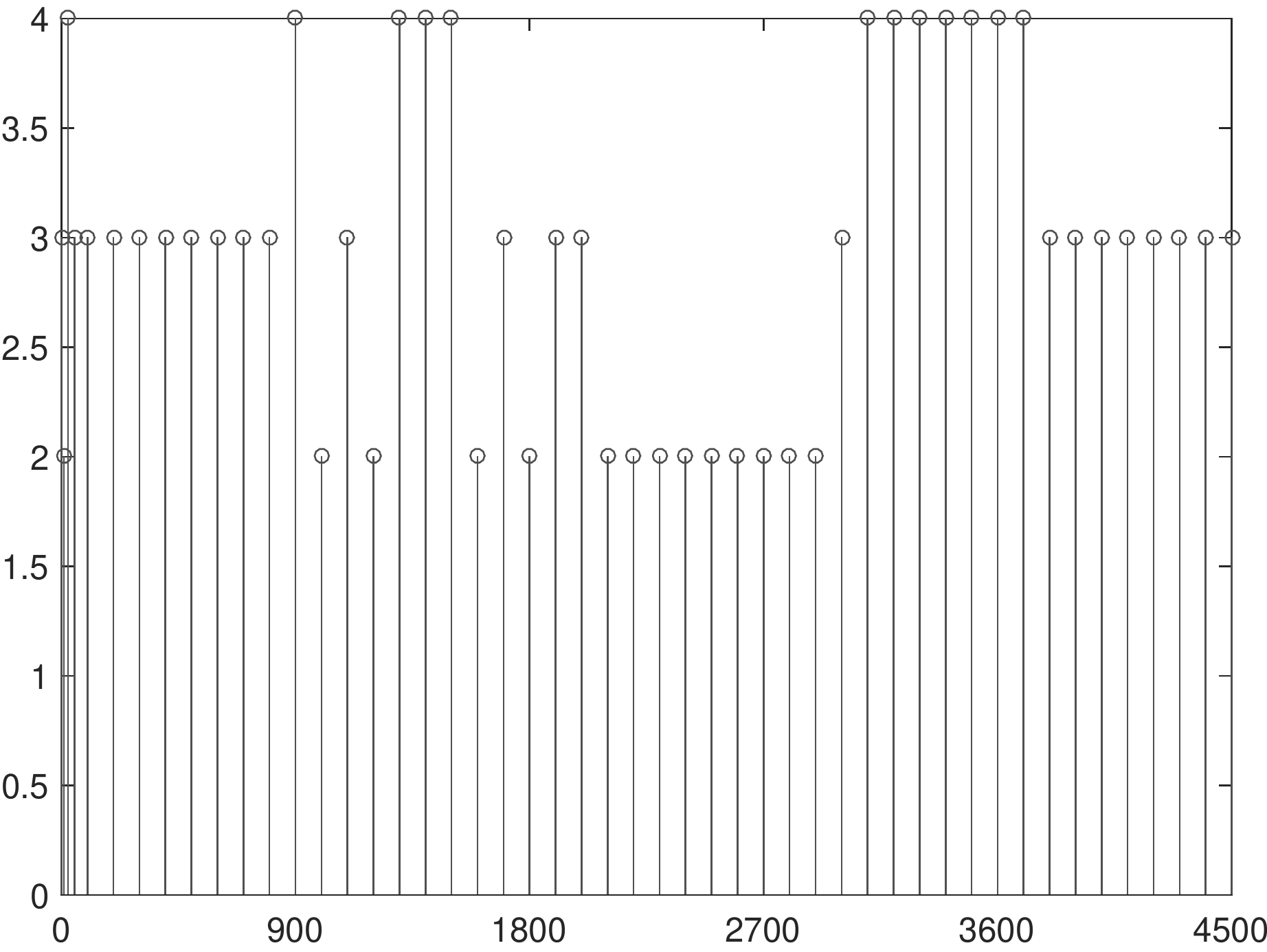}
		}
		\listsubcaptions
		\centering
		\caption{Computational Results of Algorithm \ref{alg:bernsteinRelaxedGeneralOptimizationAlgorithm}.}
		\label{fig:comp}
	\end{figure}

	All experiments are conducted on a laptop with Intel Core i7 processor with 4 physical cores and hyper-threading on each core. The maximum frequency is 2.4 GHz with the boost at specific core up to 3.2 GHz. The maximum amount of RAM allowed for computation is 2 GB for each core. Step \ref{alg:step1} in Algorithm \ref{alg:bernsteinRelaxedGeneralOptimizationAlgorithm} is coded using the optimization toolbox of Matlab R2015a, and step \ref{alg:step2} is solved via the QCP solver of CPLEX 12.6. Both of the solvers work in parallel modes, which create $ 4 $ clusters for Matlab and $ 4 $ threads for CPLEX.
	
	We run Algorithm \ref{alg:bernsteinRelaxedGeneralOptimizationAlgorithm} for case (i) with $\epsilon = 0.1$. Figure \ref{fig:comp} reports the optimal values and solutions, running times, and numbers of generated cuts when the degree $n$ of model \ref{pr:BSDP} increases from 10 to 4500. Recall that the degree of model \ref{pr:BSDP} is the degree of Bernstein polynomials approximating the set $\fU^m (\epsilon)$. Shown in Figures \ref{fig:compTotalWealth} and \ref{fig:compSol}, the optimal values and solutions fluctuate at the low degrees but become stable for the degrees larger than 3500. 
	
	The running time of Algorithm \ref{alg:bernsteinRelaxedGeneralOptimizationAlgorithm} is related to not only the degree but also the number of generated cuts. Note that the degree decides the number of decision variables in model \ref{pr:BSDP}, and the number of generated cuts is the total iterations run by Algorithm \ref{alg:bernsteinRelaxedGeneralOptimizationAlgorithm}. Figure \ref{fig:compRunningTime} reflects the tendency that a longer running time is needed as the degree increases. However, the number of generated cuts is independent of the degree shown in Figure \ref{fig:compCuts}. Particularly, Algorithm \ref{alg:bernsteinRelaxedGeneralOptimizationAlgorithm} generates 4 cuts when the degree varies in [3200, 3800], but there are only 3 cuts for the degree in [3900, 4500]. As the results, the running time reaches the peak, which is 710.11 seconds, for the degree is 3800, while it falls down to 494.94 seconds for the degree is 3900. Then the running time grows again to 647.01 seconds as the degree increases to 4500.

	\subsection{Effect Analysis of the RSD Constraint}
\begin{figure}[ht]
		\subfloat[Case (i).]{
			\label{fig:confDefault}
			\includegraphics[width=3in,height=2.2in]{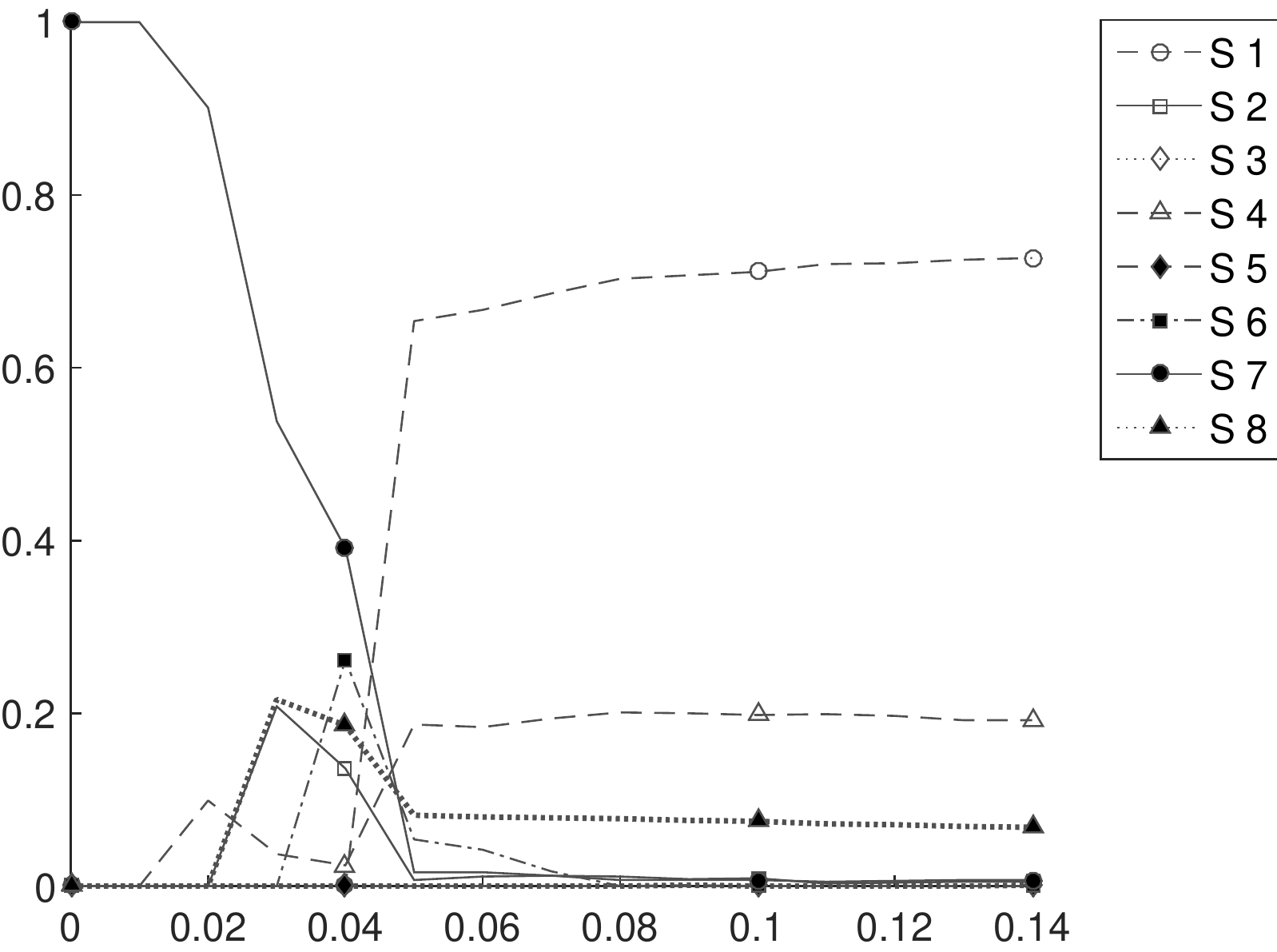}
		}
		\qquad
		\subfloat[Case (ii).]{
			\label{fig:confUref2}
			\includegraphics[width=3in,height=2.2in]{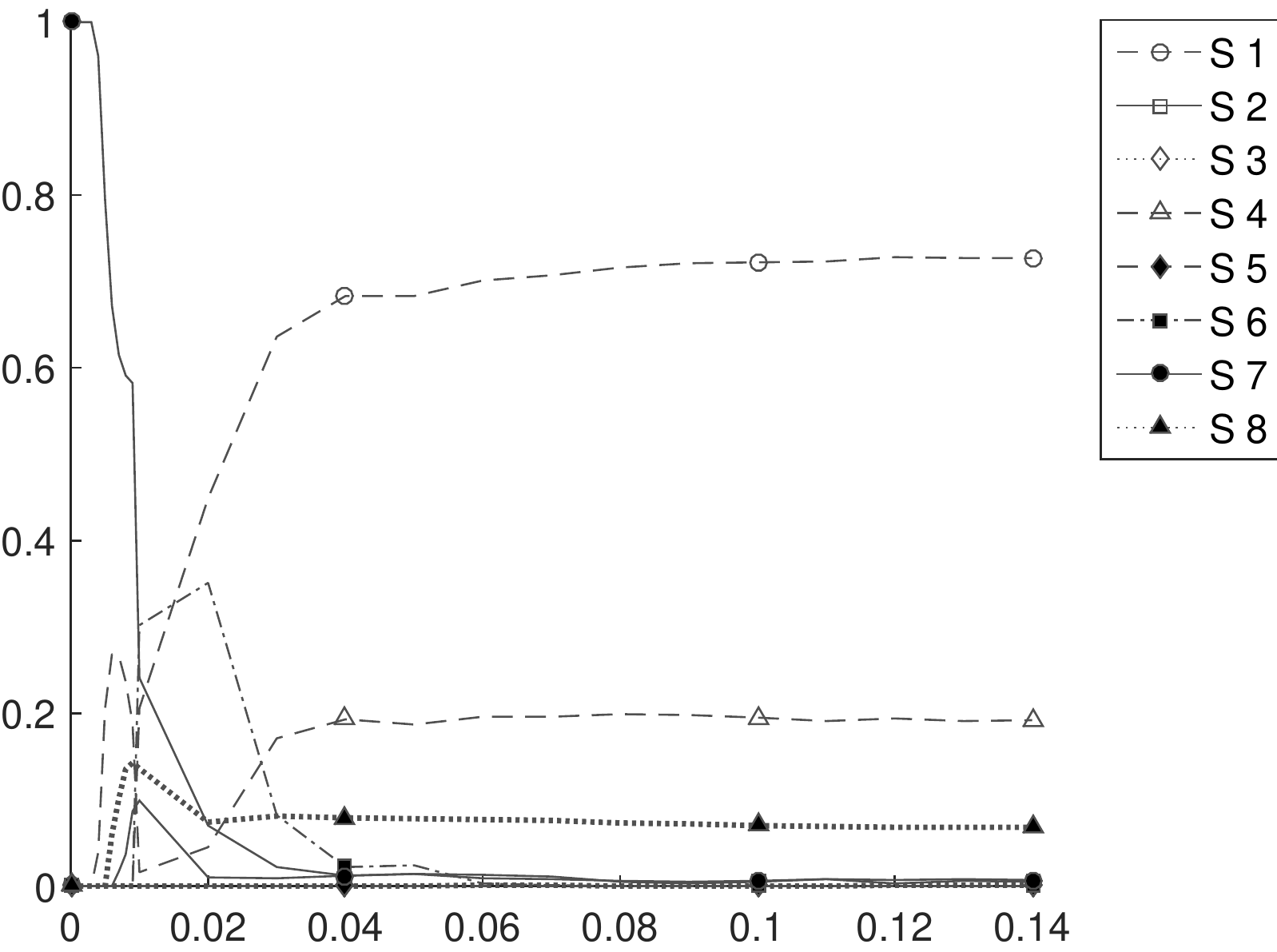}
		}
		\\[50pt]
		\subfloat[Case (iii).]{
			\label{fig:confRus}
			\includegraphics[width=3in,height=2.2in]{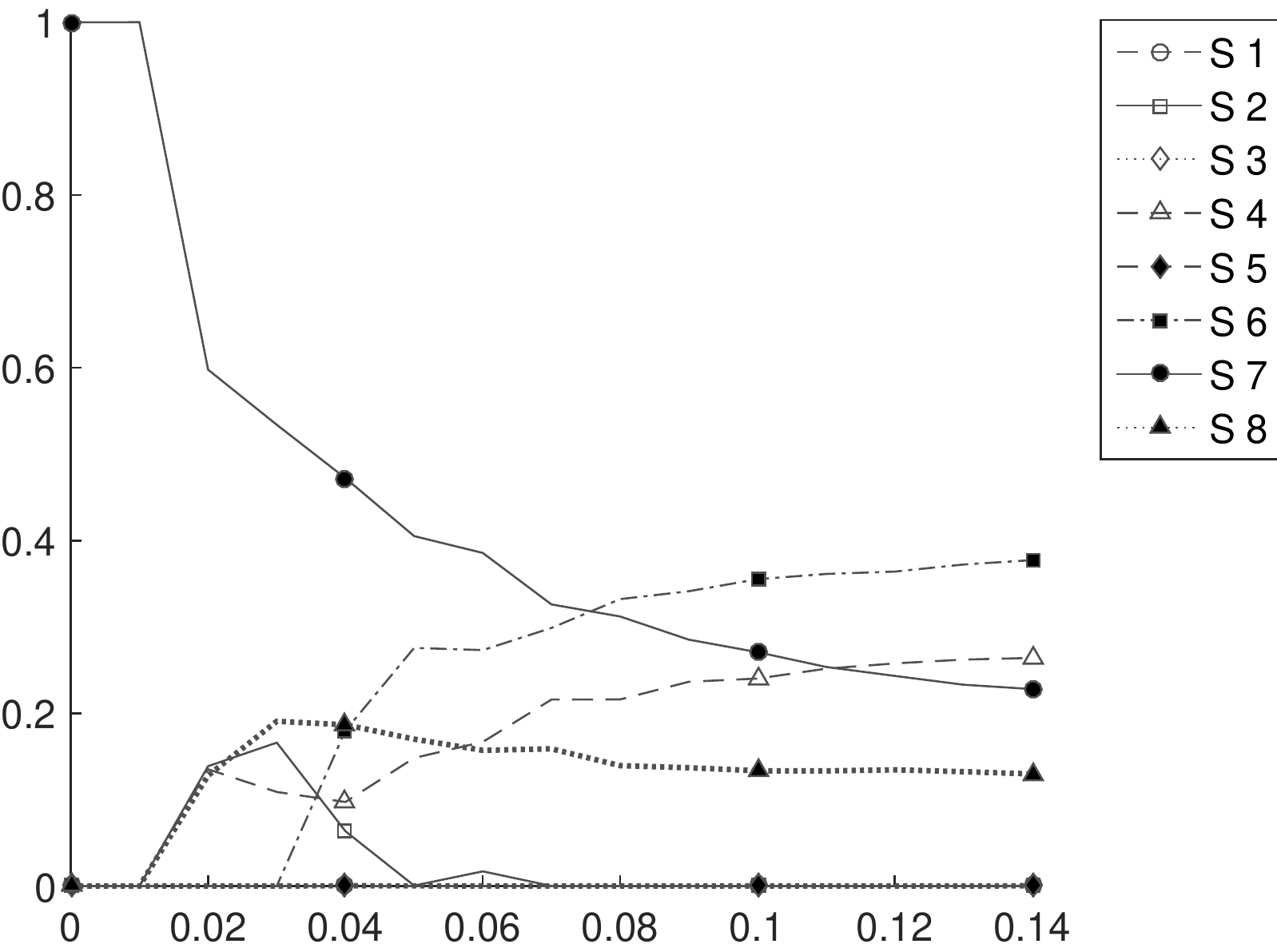}
		}
		\qquad
		\subfloat[Case (iv).]{
			\label{fig:confThird}
			\includegraphics[width=3in,height=2.2in]{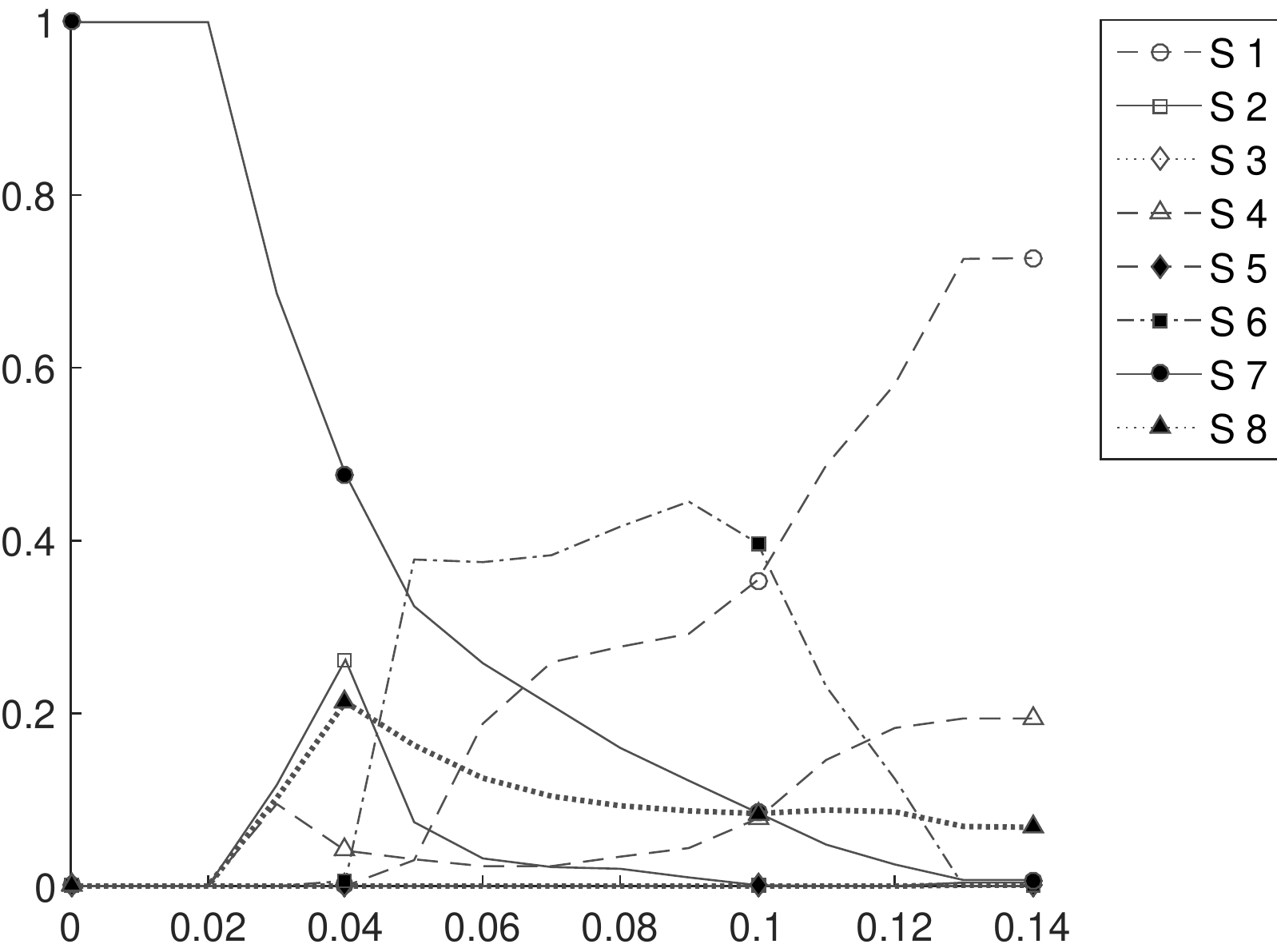}
		}
		\listsubcaptions
		\centering
		\caption{Impact of $ \epsilon $ on the Optimal Asset Allocations for Cases (i)-(iv).}
		\label{fig:config}
	\end{figure}

\begin{figure}[ht]
			\centering
			\includegraphics[width=3in]{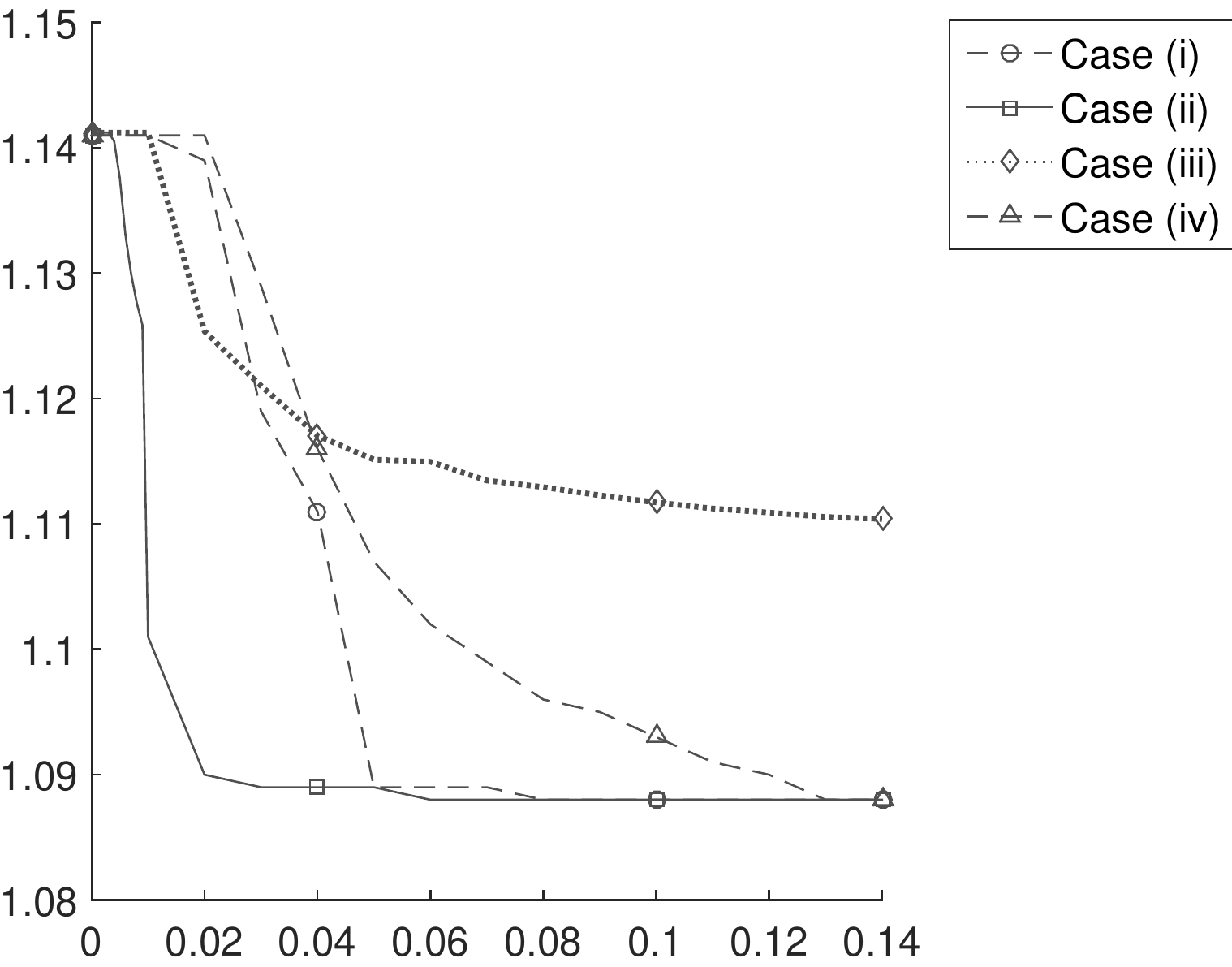}
			\caption{Expected Total Wealths of Cases (i)-(iv).}
			\label{fig:objFun}
		\end{figure}

	We now analyze cases (i) - (iv) to test the effect of the RSD constraint \eqref{pr:portfolio_RSD}. The results are given in Figures \ref{fig:config} and \ref{fig:objFun}. In this test, the degree of approximation is 4500, and
	the dominance level $\epsilon$ is adjusted in [0, 0.14]. In each case we divide this interval into three sub-intervals --- weak region, mild region, and strong region ---  due to the strength of the RSD constraint \eqref{pr:portfolio_RSD}. In general, $\epsilon$ in the weak region is very small such that the optimal value and solutions of model \eqref{pr:portfolio} are identical to ones given by only using the reference utility function (i.e., $\epsilon = 0$). Indeed, the RSD constraint \eqref{pr:portfolio_RSD} has a limited impact on the performance of model \eqref{pr:portfolio} in the weak region. The strong region is opposite, for $\epsilon$ is rather large. The corresponding optimal value and solutions are very stable, and are indifferent to the reference utility function. In contrast, model \eqref{pr:portfolio} is sensitive to $\epsilon$ in the mild region. A small change on $\epsilon$ may incur completely different asset allocations in the portfolio. 
	\\ \\
	
	\textbf{Case (i).} 
	In this case, the RSD constraint \eqref{pr:portfolio_RSD} compares risky asset allocations with the benchmark $z^{Y_1}$, which only invests the risk-free asset S1.  Shown in Figures \ref{fig:confDefault} and \ref{fig:objFun}, the weak region is [0, 0.01], the mild region is (0.01, 0.08), and the strong region is [0.08, 0.14]. Model \eqref{pr:portfolio} with $\epsilon$ in the weak region suggests that S7 obtain 100\% of the total investment and yields the highest expected total wealth 1.141. As we increase $\epsilon$ to the mild region, the investment is diversified. For example, at $\epsilon = 0.04$, the percentage of S7 in the portfolio dramatically decreases from 100\% to 39.1\%, while the percentages of S2, S4, S6, and S8 rise to 13.6\%, 2.4\%, 26.2\%, and 18.6\%, respectively. At $\epsilon = 0.06$, S1 becomes crucial in the portfolio, owning 70.3\% of the total investment and overwhelming S7 of which the percentage reduces to 1.2\%. These results reflect the fact that, for satisfying the sufficient preference over the benchmark, the RSD constraint \eqref{pr:portfolio_RSD} requires a large percentage of the total investment on the risk-free asset S1 to reduce the investment risk. As the effect of the RSD constraint \eqref{pr:portfolio_RSD} is enhanced by increasing $\epsilon$,  S1 gets more percentage until the strong region is reached. In the strong region, the investment  is stable at (72.7\%, 0.5\%, 0.1\%, 19.2\%, 0\%, 0\%, 0.7\%, 6.8\%), which is very close to the solution (72.7\%, 0.4\%, 0\%, 19.3\%, 0\%, 0\%, 0.7\%, 6.8\%)  suggested by the classical second order stochastic dominance ($\epsilon = 1$). In addition, the decrease in the investment risk greatly reduces the expected total wealth, which rapidly decreases from 1.141 to 1.089 as $\epsilon$ changes from 0 to 0.05, and then slowly changes to 1.088. The over-conservativeness of stochastic dominance results in a very low yield, in comparison with the risk-free investment on S1 yielding the expected total wealth 1.078.

	\textbf{Case (ii).} 
	This case is designed to test the effect of the reference utility function in the RSD constraint \eqref{pr:portfolio_RSD}. We substitute the CARA utility function $u^2_{ref}$ for the CRRA utility function $u^1_{ref}$. In contrast, for the total wealth more than 1, $u^2_{ref}$ has a higher Arrow and Pratt's measure of risk-aversion than $u^1_{ref}$, i.e., 
	$$
	- \frac{(u^2_{ref})''(x)}{(u^2_{ref})'(x)} =  1 > - \frac{(u^1_{ref})''(x)}{(u^1_{ref})'(x)} = \frac{1}{x}, \quad \text{for } x > 1.
	$$
	Hence, in this study, $u^2_{ref}$ characterizes stronger preference for low-risk investment than $u^1_{ref}$. It can be seen in Figures \ref{fig:confUref2} and \ref{fig:objFun} that this substitution shrinks the weak region to [0, 0.003]. A subtle change on $\epsilon = 0$ has a big impact on the investment proportion and total wealth. Analogous to case (i), the investment is diversified to hedge the risk in the mild region (0.003, 0.08). However, case (ii) has a much faster diversification rate. The asset allocation at $\epsilon = 0.04$ is (68.3\%, 1.2\%, 0\%, 19.3\%, 0\%, 2.2\%, 1.2\%, 7.9\%), in which S1 has become a major invested asset, compared to 0\% of the total investment on S1 in case (i). Also, this allocation is close to the stable solution, (72.7\%, 0.5\%, 0.1\%, 19.2\%, 0\%, 0\%, 0.7\%, 6.8\%), obtained in the strong region [0.08, 0.14]. Cases (i) and (ii) have the same asset allocation and total wealth in the strong region. This observation verifies that, for a sufficiently large $\epsilon$, the RSD constraint \eqref{pr:portfolio_RSD} is indifferent to the reference utility function, and approaches to the classical stochastic dominance.

	\textbf{Case (iii).} 
	In this case the equal allocation benchmark $z^{Y_2}$ is substituted for the risk-free investment $z^{Y_1}$. Model \eqref{pr:portfolio} suggests a completely different investment policy without over-emphasizing the safety of investment. Shown in Figure \ref{fig:confRus}, S1 is not invested on, but S7 always obtains more than 22.8\%. Figure \ref{fig:objFun} also indicates that the risky investment greatly raises the expected total wealth. In this case, the weak region is [0, 0.01] and the mild region is (0.01, 0.14]. Since the RSD constraint \eqref{pr:portfolio_RSD} is ineffective in the weak region, the best policy is still 100\% of the total investment on S7. As $\epsilon$ increases, the investment on S7 monotonically decreases, while S2, S4, S6, and S8 obtain more. The asset allocation is (0\%, 6.3\%, 0\%, 9.7\%, 0\%, 18\%, 47.2\%, 18.7\%) at $\epsilon = 0.04$, and changes to (0\%, 0\%, 0\%, 26.4\%, 0\%, 37.7\%, 22.8\%, 13\%) at $\epsilon = 0.14$. Different from case (i), there is not an overwhelming asset in case (iii), restricted by the benchmark $z^{Y_2}$ where every asset is equally treated.  
	
	\textbf{Case (iv).}
	This case tests the 3rd order RSD. As indicated by Proposition \ref{pp:orders_SD}, case (iv) with the 3rd order RSD constraint \eqref{pr:portfolio_RSD} is the relaxation of case (i). As the result, shown in Figures \ref{fig:config} and \ref{fig:objFun}, the weak region is enlarged to [0, 0.02], the strong region shrinks to [0.13, 0.14], and the diversification rate is much slower in the widely spanned mild region. Moreover, the curve of the total wealth in case (iv) is always above the curve in case (i). At $\epsilon = 0.04$, the asset allocation is (0\%, 26.2\%, 0\%, 4.1\%, 0\%, 0.6\%, 47.7\%, 21.3\%), in which case (iv) suggest 8.6\% (= 47.7\% - 39.1\%) of the total investment on S7  more than case (i). At $\epsilon = 0.06$, the allocation is (18.8\%, 3.2\%, 0\%, 2.3\%, 0\%, 37.5\%, 25.8\%, 12.5\%), in which S7 still gets more investment than S1 and S6 has the largest percentage.   
	
	\subsection{Hannah's Decision}
	We now discuss an investor's desired dominance level in this portfolio problem.  Recall that Table \ref{tab:lotteries} lists the maximum dominance level of the lottery tickets $\bar X_i$, $i = 1, \dots, 5$, with respect to non-purchase option $\bar Y_2$. We request Hannah to choose the lottery tickets which she is not reluctant to purchase, and then evaluate her desired dominance level for risky investment. Similarly, we elicit Hannah's desired dominance level for risk investment in case (i). $\bar Y_2$ does not take the risk-free interest into account. Hence,  the total wealth $W_1$ yielded by only investing on asset S1 should be substituted for $\bar Y_2$ as the risk-free investment option. Table \ref{tab:investment} gives the maximum dominance levels of lottery tickets $\bar X_i$ and $W_1$. Considering the risk-free interest leads to a little smaller levels in Table \ref{tab:investment} than in Table \ref{tab:lotteries}.
	
	\begin{table}[htbp]
		\centering
		\caption{Maximum Dominance Levels of Lottery Tickets v.s. 100\% of the Total Investment on S1}
		\begin{tabular}{crrcc}
			\addlinespace
			\toprule
			Lottery   &   \multicolumn{2}{c}{Probability of Yield}  &  & Maximum Dominance Level  \\
			Ticket & \multicolumn{1}{r}{\$0} &  \multicolumn{1}{r}{\$2} &  & $\cE^{(2)}(\bar X_{i}, W_1; u^1_{ref})$ \\
			\midrule 
			$\bar X_{1}$    &   1\%  & 99\%  &  & 0.112    \\
			$\bar X_{2}$    &   10\% & 90\%  &  & 0.061   \\
			$\bar X_{3}$    &   25\% & 85\%  &  & 0.037   \\
			$\bar X_{4}$    &   20\% & 80\%  &  & 0.018   \\
			$\bar X_{5}$    &   25\% & 75\%  &  & 0.003  \\
			\bottomrule
		\end{tabular}%
		\label{tab:investment}%
	\end{table}%
	
	Suppose that Hannah only picks up $\bar X_{1}$, and hence, her desired dominance level is 0.112. This level is in the strong region of case (i). Obviously, Hannah is a very cautious and discreet person, and thinks only $\bar X_{1}$ can almost dominate $W_1$. Figure \ref{fig:confDefault} shows that her preferred asset allocation is (72\%, 0.3\%, 0\%, 19.9\%, 0\%, 0\%, 0.5\%, 7.2\%) and the expected total wealth is 1.088. Recall that the expected total wealth is 1.078 for the risk-free investment. If both $\bar X_{1}$ and $\bar X_{2}$ can be accepted, Hannah's desired dominance level is 0.061. This level is in the mild region, and Hannah's asset allocation is (66.7\%, 1.1\%, 0\%, 18.4\%, 0\%, 4.2\%, 1.6\%, 8\%), where the investment on S1 is reduced by 5.3\% (= 72\% - 66.7\%), while the investment on  S7 increases by 1.1\% (= 1.6\% - 0.5\%). Correspondingly, the expected total wealth rises to 1.089. An appropriate assumption may be that Hannah also chooses $\bar X_{3}$ and her level is 0.037, at which the allocation is (0\%, 13.6\%, 0\%, 2.4\%, 0\%, 26.2\%, 39.1\%, 18.6\%). S7 gets 39.1\% to the total investment for ensuring a reasonable expected total wealth as 1.111, and S2, S4, S6, and S8 are chosen in the diversification for hedging risk.

\section{Conclusions}
	
	This paper has introduced a novel concept of reference-based almost stochastic dominance (RSD) and its application in risk-averse optimization problems. In the $\cL_2$-normed space, we have specified a subset of the general class of risk-averse utility functions. This subset consists of nonparametric shape-preserving perturbations around a given reference utility function. The RSD represents a preference relation that a preferred uncertain prospect should have the larger expected utility over the perturbation subset. We have also defined the maximum dominance level, which quantifies the decision maker's preference between alternative choices in the context of robustness.    
	
	We have proposed the RSD constrained stochastic optimization model and studied its solution method. An approximation approach based on Bernstein polynomials has been developed. This approach resorts to a cut-generation algorithm.	We have discussed the asymptotic convergence of the optimal value and the set of optimal solutions obtained in this approach, and proved that the algorithm has finitely many iterations. 
	
	The portfolio optimization problem given by \cite{denrus:03} has been used to analyze the computational complexity of the approximation approach and to illustrate the effect of the RSD constraint. We have compared four cases with different benchmarks, reference utility functions, and dominance orders. In addition, we have discussed the impact of an investor's desired dominance level on asset allocations.   
	
\begin{appendices}
	\section{Asset Returns given in \cite{denrus:03}}
	\label{sec:asset_returns}
	
	\begin{table}[H]
		\centering
		\caption{Asset Returns (in \%) in \cite{denrus:03}}
		\begin{tabular}{rrrrrrrrr}
			\addlinespace
			\toprule
			\multicolumn{1}{c}{Year} & \multicolumn{1}{c}{A1} & \multicolumn{1}{c}{A2} & \multicolumn{1}{c}{A3} & \multicolumn{1}{c}{A4} & \multicolumn{1}{c}{A5} & \multicolumn{1}{c}{A6} & \multicolumn{1}{c}{A7} & \multicolumn{1}{c}{A8} \\
			\midrule
			1     & 7.5   & -5.8  & -14.8 & -18.5 & -30.2 & 2.3   & -14.9 & 67.7 \\
			2     & 8.4   & 2.0   & -26.5 & -28.4 & -33.8 & 0.2   & -23.2 & 72.2 \\
			3     & 6.1   & 5.6   & 37.1  & 38.5  & 31.8  & 12.3  & 35.4  & -24.0 \\
			4     & 5.2   & 17.5  & 23.6  & 26.6  & 28.0  & 15.6  & 2.5   & -4.0 \\
			5     & 5.5   & 0.2   & -7.4  & -2.6  & 9.3   & 3.0   & 18.1  & 20.0 \\
			6     & 7.7   & -1.8  & 6.4   & 9.3   & 14.6  & 1.2   & 32.6  & 29.5 \\
			7     & 10.9  & -2.2  & 18.4  & 25.6  & 30.7  & 2.3   & 4.8   & 21.2 \\
			8     & 12.7  & -5.3  & 32.3  & 33.7  & 36.7  & 3.1   & 22.6  & 29.6 \\
			9     & 15.6  & 0.3   & -5.1  & -3.7  & -1.0  & 7.3   & -2.3  & -31.2 \\
			10    & 11.7  & 46.5  & 21.5  & 18.7  & 21.3  & 31.1  & -1.9  & 8.4 \\
			11    & 9.2   & -1.5  & 22.4  & 23.5  & 21.7  & 8.0   & 23.7  & -12.8 \\
			12    & 10.3  & 15.9  & 6.1   & 3.0   & -9.7  & 15.0  & 7.4   & -17.5 \\
			13    & 8.0   & 36.6  & 31.6  & 32.6  & 33.3  & 21.3  & 56.2  & 0.6 \\
			14    & 6.3   & 30.9  & 18.6  & 16.1  & 8.6   & 15.6  & 69.4  & 21.6 \\
			15    & 6.1   & -7.5  & 5.2   & 2.3   & -4.1  & 2.3   & 24.6  & 24.4 \\
			16    & 7.1   & 8.6   & 16.5  & 17.9  & 16.5  & 7.6   & 28.3  & -13.9 \\
			17    & 8.7   & 21.2  & 31.6  & 29.2  & 20.4  & 14.2  & 10.5  & -2.3 \\
			18    & 8.0   & 5.4   & -3.2  & -6.2  & -17.0 & 8.3   & -23.4 & -7.8 \\
			19    & 5.7   & 19.3  & 30.4  & 34.2  & 59.4  & 16.1  & 12.1  & -4.2 \\
			20    & 3.6   & 7.9   & 7.6   & 9.0   & 17.4  & 7.6   & -12.2 & -7.4 \\
			21    & 3.1   & 21.7  & 10.0  & 11.3  & 16.2  & 11.0  & 32.6  & 14.6 \\
			22    & 4.5   & -11.1 & 1.2   & -0.1  & -3.2  & -3.5  & 7.8   & -1.0 \\
			\bottomrule
		\end{tabular}%
		\label{tbl:Index_Return}%
	\end{table}%
\end{appendices}

	
\bibliographystyle{ormsv080} 
\bibliography{draft}

\end{document}